%% file: RR-6595.tex
\documentclass[a4paper,11pt]{article}
\usepackage{color}
\usepackage{float}
\usepackage{RR}
\usepackage{hyperref}
\usepackage{pstricks,pst-node,pst-text,pst-3d}
\usepackage{rotate,float,pst-plot}
\usepackage{latexsym}
\usepackage{amssymb}
\usepackage[latin1]{inputenc}
\usepackage{fancyhdr}
\usepackage{amsmath}
\usepackage{amsthm}
\usepackage{a4wide}
\usepackage[amssymb,Gray]{SIunits}
\RRdate{Juillet 2008} 

\RRauthor{
Julien Diaz
  \thanks[sfn]{EPI Magique-3D, Centre de Recherche
  Inria Bordeaux Sud-Ouest}
\thanks[sf1]{Laboratoire de Mathématiques et de
  leurs Applications, CNRS UMR-5142, Université de Pau et des Pays de
  l'Adour --  Bâtiment IPRA, avenue de
  l'Université  -- BP 1155-64013 PAU CEDEX}%
  \and
Abdelaâziz Ezziani
\thanksref{sf1}
\thanksref{sfn}
}
\authorhead{Diaz \& Ezziani}
\RRtitle{Solution analytique pour la propagation
  d'ondes en milieu stratifié
hétérogène acoustique/poroélastique. Partie~II~: en dimension 3}
\RRetitle{Analytical Solution for Wave Propagation in Stratified
Acoustic/Porous Media. Part~II: the 3D Case} 
\titlehead{Analytical solution}
\RRresume{Nous nous intéressons à la modélisation de la propagation
d'ondes dans les milieux infinis bicouche
acoustique/poroélastique.  Nous considérons le modèle
bi-phasique de Biot dans la couche
poroélastique. La première partie est
consacrée au calcul de la solution analytique en
dimension deux à l'aide de la
technique de Cagniard-De Hoop. Dans cette deuxième partie nous considérons le 
cas de la dimension 3.
}
  \RRabstract{We are interested in the modeling of
    wave propagation in an infinite bilayered acoustic/poroelastic media.
  We consider the biphasic Biot's model in the
  poroelastic layer. The first part is devoted to the
  calculation of analytical solution in two
  dimensions, thanks to  Cagniard de Hoop method. In this second part 
we consider the 3D case.
} 
\RRmotcle{Modèle de Biot, ondes poroélastiques,
  ondes acoustiques, couplage
  acoustique/poroelastique, solution analytique,
technique de Cagniard de Hoop.}
 \RRkeyword{Biot's model, poroelastic waves,
   acoustic waves, acoustic/poroelastic coupling, analytical solution, Cagniard-De Hoop's technique.} 
 \RRprojet{Magique-3D}  
\RRtheme{\THNum} 
 \URFuturs 
\RCBordeaux
\input{ncacousporo}
\RRNo{6595}
\begin{document}
\makeRR
\section*{Introduction}
The computation of analytical solutions for wave
propagation problems is of high importance for the
validation of numerical computational codes or for
a better understanding of the
reflexion/transmission properties of the media.
Cagniard-de Hoop method~\cite{Cag,DH} is a useful
tool to obtain such solutions and permits to
compute each type of waves (P wave, S wave, head
wave...) independently.  Although it was
originally dedicated to the solution of
elastodynamic wave propagation, it can be applied
to any transient wave propagation problem in
stratified media. However, as far as we know,
 few works have been dedicated to the
application of this method to poroelastic medium.
In~\cite{ezziani_th} the analytical solution of
poroelastic wave propagation in an homogeneous 2D
medium is provided and in~\cite{FJ} the authors
compute the analytical expression of the reflected
wave at the interface between an acoustic and a
poroelastic layer in two dimension but they do not
explicit the expression of the transmitted waves.

In order to validate computational codes of wave
propagation in poroelastic media, we have
implemented the codes Gar6more 2D~\cite{Gar6} and
Gar6more 3D~\cite{Gar63d} which provide the complete solution
(reflected and transmitted waves) of the
propagation of wave in stratified 2D or 3D media
composed of acoustic/acoustic, acoustic/elastic,
acoustic/poroelastic or poroelastic/poroelastic
The 2D code and the 3D code are freely downloadable at\\
\centerline{\url{http://www.spice-rtn.org/library/software/Gar6more2D}.}
and\\
\centerline{\url{http://www.spice-rtn.org/library/software/Gar6more3D}.}
In previous studies~\cite{RAP_DE6509,RAP3} we have
presented the 2D acoustic/poroelastic and poroelastic/poroelastic cases and  we  focus here
on the 3D acoustic/poroelastic case,the
3D poroelastic/poroelastic case will be
the object of forthcoming papers.  We first present the
model problem we want to solve and derive the
Green problem from it (section 1).  Then we
present the analytical solution of wave
propagation in a stratified 3D medium composed of
an acoustic and a poroelastic layer (section 2)
and we detail the computation of the solution
(section 3). Finally we illustrate our results
through numerical applications (section 4).

\section{The model problem}
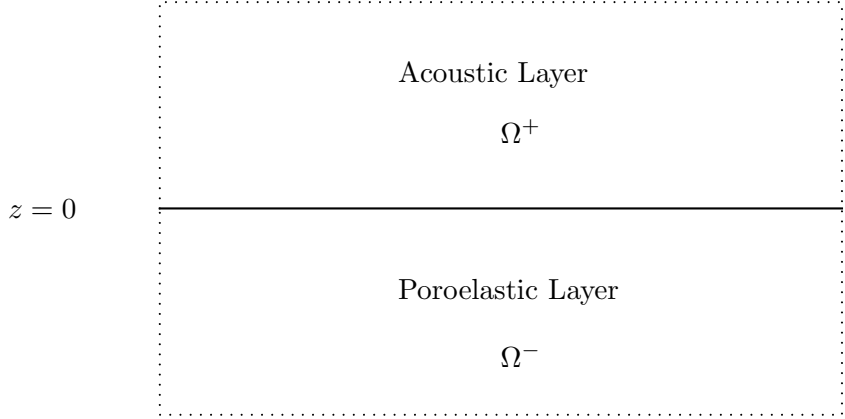
\begin{figure}[htbp]
  \centering
\setlength{\unitlength}{.9cm}
  \begin{picture}(6,6.5)(2,0)
\psset{xunit=1.5cm}
\psframe[linestyle=dotted](0,.25)(6,5.75)
\psline(0,3)(6,3)
\put(5,4.3){$\Omega^+$}
\put(5,1.){$\Omega^-$}
\put(3.5,5.2){Acoustic Layer}
\put(3.5,2){Poroelastic Layer}
\put(-2.2,3.2){$z=0$}
  \end{picture}
\caption{Configuration of the study}
\label{fig:interf_plan}
\end{figure}

We consider an infinite three dimensional medium
($\Om=\R^3$) composed of an homogeneous acoustic layer
$\Omega^+=\R^3\times]-\infty,0]$ and an homogeneous poroelastic
layer $\Omega^-=\R^3\times[0,+\infty[$ separated
by an horizontal interface $\Gamma$ (see
Fig.~\ref{fig:interf_plan}). We first describe the
equations in the two layers
(\S\ref{sec:equation-acoustics} and
\S\ref{sec:biots-model} ) and the transmission
conditions on the interface $\Gamma$
(\S\ref{sec:transm-cond}), then we present the
Green problem from which we compute the
analytical solution (\S\ref{sec:greens-problem}).
\subsection{The equation of acoustics}
\label{sec:equation-acoustics}
In the acoustic layer we consider the second order
formulation of the wave equation with a point
source in space, a regular source function $f$
in time and zero initial conditions:
\begin{equation}
\label{eq:acous}
\left\{
\begin{array}{ll}
\ddot P^+ -{V^+}^2 \Delta
P^+=\delta_x\,\delta_y\,\delta_{z-h}\,f(t),&\mbox{in }
\Om^+\times]0,T],\\[8pt]
\dsp \ddot {\gr U}^+=-\frac{1}{\rho^+}\nabla P^+,&\mbox{in }
\Om^+\times]0,T],\\[8pt]
P^+(x,y,0)=0,\dot P^+(x,y,0)=0 ,&\mbox{in }
\Om^+\\[8pt]
{\bf U}^+(x,y,0)=0,\dot {\bf U}^+(x,y,0)=0 ,&\mbox{in }
\Om^+
\end{array}
\right.
\end{equation}
where
\begin{itemize}
\item $P^+$ is the pressure;\\
\item $\gr U^+$ is the displacement field;\\
\item $V^+$ is the celerity of the wave;\\
\item $\rho^+$ is the density of the fluid.
\end{itemize}
\subsection{Biot's model}
\label{sec:biots-model}
In the second layer we consider the second order
formulation of the poroelastic
equations~\cite{biot1,biot2,biot3}
\begin{equation}\label{eq:biot}
\left\{
\begin{array}{lll}
\dsp\ro^-\,\ddU^-+\ro_f^-\,\ddW^--\Nabla\cdot\Si^-=0,&&\mbox{in }
\Om^-\times]0,T],\\[8pt]
\dsp\ro_f^-\,\ddU^-+\ro_w^-\,\ddW^-+\frac{1}{\mat K^-}\,\gr \dW^-+\nabla P^-=0,
&&\mbox{in }\Om^-\times]0,T],\\[8pt]
\dsp \Si^-=\la^-\nabla\cdot \gr U^-_s\,\gr I_3+2\mu^-\vare(\gr U^-_s)-\be^-\, P^-\,\gr I_3, &&\mbox{in }\Om^-\times]0,T],
\\[8pt]
\dsp\frac{1}{m^-}\,P^-+\be^-\,\nabla\cdot\gr U^-_s+\nabla\cdot \gr
W^-=0,&&\mbox{in }\Om^-\times]0,T],\\[12pt]
\gr U^-_s(x,0)=0,\,\gr W^-(x,0)=0, &&
\mbox{in }\Om^-,
\\[12pt]
\dU^-(x,0)=0,\,\dW^-(x,0)=0, &&
\mbox{in }\Om^-,
\end{array}
\right.
\end{equation}
with
$$
(\Nabla\cdot\Si^-)_i=\sum_{j=1}^3\frac{\p\Si^-_{ij}}{\p
  x_j}\;\;\forall\,i=1,3,\;\gr{I}_3 \mbox{ is the
  usual identity matrix of }\mat M_3(\RR),
$$
and $\vare(\gr U^-_s)$ is the solid strain tensor defined by: 
$$
\vare_{ij}(\gr U)=\fr\left(\frac{\p U_i}{\p x_j}+\frac{\p U_j}{\p x_i}\right).$$
In (\ref{eq:biot}), the unknowns are:\\
\begin{itemize}
\item $\gr U^-_s$ the displacement field of solid particles;\\
\item $\gr W^-=\phi^-(\gr U^-_f-\gr U^-_s)$, the relative displacement, $\gr
  U^-_f$ being the displacement field of fluid particle and $\phi^-$ 
the porosity;\\
\item $P^-$, the fluid pressure;\\
\item $\Si^-$, the solid stress tensor.\\
\end{itemize}
The parameters describing the physical properties
of the medium are 
given by:\\
\begin{itemize}
\item $\ro^-=\phi^-\,\ro^-_f+(1-\phi^-)\ro^-_s$ is the overall density of the saturated medium, with $\ro^-_s$ the density of the solid and  $\ro^-_f$ the density of the fluid;\\
\item $\ro_w^-=a^-\ro^-_f/\phi^-$, where $a^-$ the tortuosity of the solid matrix;\\
\item $\mat K^-=\kappa^-/\eta^-$, $\kappa^-$ is the permeability of the solid matrix and $\eta\-$ is the viscosity of the fluid;\\
\item $m^-$ and $\beta^-$ are positive physical coefficients: 
$\be^-=1-K^-_b/K^-_s$ \\
and $m^-=\left[\phi^-/K^-_f+(\be^--\phi^-)/K_s^-\right]^{-1}$, where 
$K^-_s$ is the bulk modulus of the solid, $K^-_f$ is the bulk modulus of the fluid 
and $K^-_b$ is the frame bulk modulus;\\
\item  $\mu^-$ is the frame shear modulus,
and $\la^-=K_b^--2\mu^-/3$ is the Lam\'e constant.
\end{itemize}
\subsection{Transmission conditions}
\label{sec:transm-cond}
Let $\gr n$ be the unitary normal vector of
$\Gamma$ outwardly directed to. The transmission 
conditions on the interface between the acoustic and porous medium  
are~\cite{carcione} :
\begin{equation}\label{eq:condi_trans}
\left\{
\begin{array}{l}
\gr W^- \cdot \gr n=(\gr U^+-\gr U_s^-)\cdot\gr n,\\[4pt]
P^-=P^+,\\[4pt]
\Si^-\,\gr n=-P^+\gr n.
\end{array}
\right.
\end{equation}
\subsection{The Green problem}
\label{sec:greens-problem}
We won't compute directly the solution to
(\ref{eq:acous}-\ref{eq:biot}-\ref{eq:condi_trans})
but the solution to the following Green problem:
\begin{subequations}\label{eq:acous2}
\begin{eqnarray}
\ddot p^+ -{V^+}^2 \Delta
p^+=\delta_x\,\delta_y\,\delta_{z-h}\,\delta_t,&&\mbox{in }
\Om^+\times]0,T],\label{eq:acous1}\\[8pt]
\dsp \ddot {\gr u}^+=-\frac{1}{\rho^+}\nabla p^+,&&\mbox{in }
\Om^+\times]0,T],\label{eq2:acous1}\
\end{eqnarray}
\end{subequations}
\begin{subequations}\label{biot_decomp}
\begin{eqnarray}
\dsp\ro^-\,\ddu^-+\ro_f^-\,\ddw^--\Nabla\cdot\si^-=0,&&\mbox{in }
\Om^-\times]0,T],\label{eq:biot1}\\[8pt]
\dsp\ro_f^-\,\ddu^-+\ro_w^-\,\ddw^-+\frac{1}{\mat K^-}\,\gr \dw^-+\nabla p^-=0,
&&\mbox{in }\Om^-\times]0,T],\label{eq:biot2}\\[8pt]
\dsp \si^-=\la^-\nabla\cdot \gr u^-_s\,\gr I_3+2\mu^-\vare(\gr u^-_s)- \be^-\, p^-\,\gr I_3, &&\mbox{in }\Om^-\times]0,T],
\label{acousporo:eq:6}\\[8pt]
\dsp\frac{1}{m^-}\,p^-+\be^-\,\nabla\cdot\gr u^-_s+\nabla\cdot \gr
w^-=0,&&\mbox{in }\Om^-\times]0,T],\label{acousporo:eq:8}
\end{eqnarray}
\end{subequations}
\begin{subequations}\label{eq:condi_trans2}
\begin{eqnarray}
\gr w^-\cdot\gr n=(\gr u^+-\gr u_s^-)\cdot\gr n,&&\mbox{on }\Gamma,\\[4pt]
p^-=p^+,&&\mbox{on }\Gamma,\\[4pt]
\si^-\,\gr n=-p^+\,\gr n,&&\mbox{on }\Gamma.
\end{eqnarray}
\end{subequations}
The solution to (\ref{eq:acous}-\ref{eq:biot}-\ref{eq:condi_trans})
is then computed from the solution of the Green
Problem thanks to a  convolution by the source function. For instance we have~:
$$P^+(x,y,t)=p^+(x,y,.)\ast f(.)=\int_0^t p^+(x,y,\tau)f(t-\tau)\,d\tau$$
(we have similar relations for the other
unknowns). We also suppose that the poroelastic
medium is non dissipative, i.e the viscosity
$\eta^-=0$. Using the equations
(\ref{acousporo:eq:6}, \ref{acousporo:eq:8}) we
can eliminate $\si^-$ and $p^-$ in
(\ref{biot_decomp}) and we obtain the the equivalent system:  
\begin{equation}\label{syst-biot2}
\hspace*{-0.2cm}\left\{
\begin{array}{ll}
\ro^-\,\ddu^-+\ro_f^-\,
\ddw^--\alpha^-\,\nabla(\nabla\cdot \gr u_s^-)
+\mu^-\,\nabla\times(\nabla\times \gr u_s^-)-m^-\be^-\nabla(\nabla\cdot \gr w^-)=0,&z<0\\[12pt]
\ro_f^-\,\ddus^-
+\ro_w^-\,\ddw^--m^-\be^-\,\nabla(\nabla\cdot\gr u_s^-)
-m^-\,\nabla(\nabla\cdot\gr w^-)=0,&z<0
\end{array}
\right .
\end{equation}
with $\alpha^-=\la^-+2\mu^-+m^-{\be^-}^2$.\\\\
And using the equation (\ref{eq2:acous1}) the transmission conditions 
 (\ref{eq:condi_trans2}) on $z=0$ are rewritten as:
\begin{subequations}\label{eq:condi_trans_expli}
\begin{eqnarray}
&&\dsp \ddot u_{s\,z}^-+\ddot w_{z}^-=-\frac{1}{\ro^+}\p_z p^+,
\label{eq:condi_trans_expli1}\\[12pt]
&&-m^-\be^-\,\nabla\cdot\gr u_s^--m^-\,\nabla\cdot\gr w^-=p^+,
\label{eq:condi_trans_expli2}\\[12pt]
&&\dsp\p_zu_{sx}^-+\p_x u_{sz}^-=0,
\label{eq:condi_trans_expli3}\\[12pt]
&&\dsp\p_zu_{sy}^-+\p_y u_{sz}^-=0,
\label{eq:condi_trans_expli4}\\[12pt]
&&\dsp(\la^-+m^-{\be^-}^2)\nabla\cdot\gr u_s^-+2\mu^-\p_z u_{sz}^-
+m^-\be^-\,\nabla\cdot\gr w=-p^+\label{eq:condi_trans_expli5}.
\end{eqnarray}
\end{subequations}
We split the displacement fields $\gr u_s^-$ and
$\gr u_f^-$ into irrotationnal and isovolumic fields
(P-wave and S-wave):
\begin{equation}\label{eq:irriso}
\gr u_s^-=\nabla \Theta_u^-+\nabla\times\gr\Psi_u^-
\;\; ; \;\;\gr w^-= \nabla \Theta_w^-+\nabla\times\gr\Psi_w^-.
\end{equation}
The vectors $\gr \Psi_u$ and $\gr \Psi_w$ are not uniquely defined since:
$$
\nabla\times (\gr \Psi_\ell+\nabla C)=\nabla\times\gr\Psi_\ell,\;\;\;\forall\,\ell\in\{u,w\}
$$  
for all scalar field $C$. To define a unique $\gr \Psi_\ell$ we impose the gauge condition:
$$
\nabla\cdot\gr \Psi_\ell=0
$$  
The vectorial space of $\gr \Psi_\ell^\pm$ verifying this last condition is written as:
$$
\gr \Psi_\ell=\left[\begin{array}{c}
\p_y\\[5pt]
-\p_x\\[5pt]
0\end{array}\right]\Psi_{\ell,1}+
\left[\begin{array}{c}
\p_{xz}^2\\[5pt]
\p_{yz}^2\\[5pt]
-\p_{xx}^2-\p_{yy}^2\end{array}\right]\Psi_{\ell,2},
$$ 
where $\Psi_{\ell,1}$ and $\Psi_{\ell,2}$ are two scalar fields. The displacement fields $\gr u_s^-$ and $\gr w^-$ are written in the form:
\begin{equation}
\begin{array}{l}
\gr u_s^-=\nabla\Theta_u^-+\left[\begin{array}{c}
\p_{xz}^2\\[5pt]
\p_{yz}^2\\[5pt]
-\p_{xx}^2-\p_{yy}^2\end{array}\right]\Psi_{u,1}^--\left[\begin{array}{c}
\p_y\\[5pt]
-\p_x\\[5pt]
0\end{array}\right]\Delta\Psi_{u,2}^-\\[35pt]
\gr w^-=\nabla\Theta_w^-+\left[\begin{array}{c}
\p_{xz}^2\\[5pt]
\p_{yz}^2\\[5pt]
-\p_{xx}^2-\p_{yy}^2\end{array}\right]\Psi_{w,1}^--\left[\begin{array}{c}
\p_y\\[5pt]
-\p_x\\[5pt]
0\end{array}\right]\Delta\Psi_{w,2}^-
\end{array}
\end{equation}
We can then rewrite system (\ref{syst-biot2}) in the following form:
\begin{equation}\label{equ:matrice}
\left\{
\begin{array}{ll}
A^-\ddot\Theta^--B^-\Delta \Theta^-=0,&z<0\\[8pt]
\ddot\Psi_{u,1}^--{\Vs^-}^2\Delta \Psi_{u,1}^-=0,&z<0\\[8pt]
\ddot\Psi_{u,2}^--{\Vs^-}^2\Delta
\Psi_{u,2}^-=0,&z<0\\[8pt]
\dsp
\ddot{\gr\Psi}_w^-=-\frac{\ro_f^-}{\ro_w^-}\ddot{\gr\Psi}_u^-,&z<0\end{array}
\right.
\end{equation}
where $\Theta^-=(\Theta_u^-,\Theta_w^-)^t$, $A^-$
and $B^-$ are $2\times 2$ symmetric matrices:
$$
A^-=\left(\begin{array}{cc}
\ro^-&\ro_f^-\\[8pt]
\ro_f^-&\ro_w^-
\end{array}\right)\;\;;\;\;B^-=\left(\begin{array}{cc}
\lambda^-+2\mu^-+m^-(\beta^-)^2&m^-\beta^-\\[8pt]
m^-\beta^-&m^-
\end{array}\right),
$$
and $$
\Vs^-=\sqrt{\frac{\mu\ro_w^-}{\ro^-\ro_w^--{\ro_f^-}^2}}$$
is the S-wave velocity.\\\\
We multiply the first equation of the system (\ref{equ:matrice}) by
the inverse of $A$. The matrix  $A^{-1}B$ is diagonalizable:
$A^{-1}B=\mat P D{\mat P}^{-1}$, where $\mat
P$ is the change-of-coordinate matrix,
$D=diag({\Vpf^-}^2,{\Vps^-}^2)$ is the diagonal matrix similar
to $A^{-1}B$, $\Vpf^-$ and $\Vps^-$ are
respectively the fast P-wave velocity and the slow
P-wave velocity ($\Vps<\Vpf$).\\\\
Using the change of variables 
\begin{equation}\label{changethetaphi}
\Phi^-=(\Phi_{Pf}^-,\Phi_{Ps}^-)^t={\mat P}^{-1}\Theta^-,
\end{equation} we obtain the uncoupled system on
fast P-waves, slow P-waves and S-waves:
\begin{equation}\label{equ:matricediag}
\left\{
\begin{array}{ll}
\ddot\Phi^--D\Delta \Phi^-=0,&z<0\\[8pt]
\ddot\Psi_{u,i}^--{\Vs^-}^2\Delta \Psi_{u,i}^-=0,\ \ i=1,2,&z<0\\[8pt]
\dsp
\Psi_w^-=-\frac{\ro_f^-}{\ro_w^-}\Psi_u^-.&z<0
\end{array}
\right.
\end{equation}
Using the transmission conditions (\ref{eq:condi_trans_expli3})-(\ref{eq:condi_trans_expli4}), we obtain:
\begin{subequations}\label{eq:trans2eq}
\begin{eqnarray}
2\p_{xz}^2\Theta_u^-+\p_x(\p_{zz}^2-\Delta_\perp)\Psi_{u,1}^-
-\p_{yz}^2\Delta\Psi_{u,2}^-=0\label{eq:trans2eq1},&&\mbox{on }\Ga,\\[8pt]
2\p_{yz}^2\Theta_u^-+\p_y(\p_{zz}^2-\Delta_\perp)\Psi_{u,1}^-
+\p_{xz}^2\Delta\Psi_{u,2}^-=0\label{eq:trans2eq2},&&\mbox{on }\Ga,
\end{eqnarray}
\end{subequations}
with $\Delta_\perp=\p_{xx}^2+\p_{yy}^2$. Applying the  derivative $\p_y$ 
to the equation (\ref{eq:trans2eq1}), $\p_x$ to
the equation (\ref{eq:trans2eq2}) and 
subtracting the first obtained equation from the second one, we get:
\begin{equation}\label{eq:psi2}
(\p_z\Delta_\perp)\Delta \Psi_{u,2}^-=0, \mbox{ on }\Ga,
\end{equation}
moreover, using that $\Psi_{u,2}^-$ satisfies the wave equation:
$$
\ddot\Psi_{u,2}^--{\Vs^-}^2\Delta \Psi_{u,2}^-=0, \;\;z<0
$$
and that $\gr u_s^-$ and $\gr w^-$ satisfy, at  $t=0$, 
$
\gr u_s^-=\dot{\gr u}_s^-=\gr w^-=\dot{\gr w}^-=0,
$
we obtain: 
$$
\Psi_{u,2}^-=0,\;z\le 0, 
$$ 
and from (\ref{eq:trans2eq}) we deduce the transmission condition equivalent to (\ref{eq:condi_trans_expli3}) and(\ref{eq:condi_trans_expli4}):
\begin{equation}\label{condtransreduite}
2\p_{z}^2\Theta_u^-+(\p_{zz}^2-\Delta_\perp)\Psi_{u,1}^-=0,\;\;\mbox{on }\Ga.
\end{equation}
Finally, we obtain the Green problem equivalent to (\ref{eq:acous2},\ref{biot_decomp},\ref{eq:condi_trans2}):
\begin{equation}\label{equ:diagsyst}
\left\{
\begin{array}{ll}
\ddot p^+ -{V^+}^2 \Delta
p^+=\delta_x\,\delta_{y}\,\delta_{z-h}\,\delta_t,&z>0\\[8pt]
\ddot\Phi_i^--{V_i^-}^2\Delta \Phi_i^-=0,\quad i\in\{Pf,Ps,S\}&z<0\\[8pt]
\dsp {\cal B} (p^+,\Phi_{Pf}^-,\Phi_{Ps}^-,\Phi_S^-)=0,&z=0
\end{array}
\right.
\end{equation}
where we have set $\Phi_{S}^-=\Psi_{u,1}^-$ in order to have similar notations for the $Pf$, $Ps$ and $S$ waves.  The operator ${\cal B}$ represents the
transmission conditions on $\Gamma$:
$$
{\cal B} \left(
\begin{array}{l}
p^+\\
\Phi_{Pf}^-\\
\Phi_{Ps}^-\\
\Phi_{S}^-\\
\end{array}\right)=\left[
\begin{array}{cccc}
\dsp \frac{1}{\ro^+}\p_z&(\mat P_{11}+\mat P_{21})\,\p_{ztt}^3
&(\mat P_{12}+\mat P_{22})\,\p_{ztt}^3&\dsp (\frac{\ro_f^-}{\ro_w^-}-1)\,
\p_{tt}^2\Delta_\perp\\[18pt]
1&\dsp \frac{m^-(\be^-\mat P_{11}+\mat P_{21})}{{\Vpf^-}^2}\p_{tt}^2&
\dsp \frac{m^-(\be^-\mat P_{12}+\mat P_{22})}{{\Vps^-}^2}\p_{tt}^2&0\\[18pt]
0&2\mat P_{11}\,\p_{z}&2\mat P_{12}\,\p_{z}&\p_{zz}^2-\Delta_\perp\\[10pt]
\dsp 1&{\cal B}_{42}&{\cal B}_{43}&-2\mu^-\p_{z}\Delta_\perp
\end{array}
\right]\left[
\begin{array}{l}
p^+\\[18pt]
\Phi_{Pf}^-\\[18pt]
\Phi_{Ps}^-\\[10pt]
\Phi_{S}^-
\end{array}
\right]
$$
where $\mat P_{ij}$, $i,j=1,2$ are the components of the change-of-coordinates matrix $\mat P$, $\mat B_{42}$ and $\mat B_{43}$ are given by: 
$$
\begin{array}{l}
{\cal B}_{42}=\dsp \frac{(\la^-+m^-{\be^-}^2)\mat P_{11}+
m^-\be^-\mat P_{21}}{{\Vpf^-}^2}\p_{tt}^2
+2\mu^-\mat P_{11}\p_{zz}^2,\\[16pt]
{\cal B}_{43}=\dsp \frac{(\la^-+m^-{\be^-}^2)\mat P_{12}+
m^-\be^-\mat P_{22}}{{\Vps^-}^2}\,\p_{tt}^2+2\mu^-\mat P_{12}\p_{zz}^2.
\end{array}
$$
To obtain this operator we have used the
transmission conditions
(\ref{eq:condi_trans_expli1},\ref{eq:condi_trans_expli2},\ref{condtransreduite},\ref{eq:condi_trans_expli5}),
the change of variables (\ref{eq:irriso}) and the
uncoupled system (\ref{equ:matricediag}).\\\\
Moreover, we can determine the solid displacement $\gr u_s^-$ by using the change of variables (\ref{eq:irriso}) and the fluid displacement $\gr u^+$ by using 
(\ref{eq2:acous1}).
\section{Expression of the analytical solution}
Since the problem is invariant by a rotation
around the $z$-axis, we will only consider the
case $y=0$ and $x>0$, so that the y-component of all the
displacements are zero. The solution for $y\neq0$
or $x\leq0$
is deduced from the solution for $y=0$ by the
relations
\begin{eqnarray}
  \label{acousporo3d:eq:2}
&&\dsp   p(x,y,z,t)=p(\sqrt{x^2+y^2},0,z,t)\\[10pt]
&&\dsp
u_{s\,x}(x,y,z,t)=\frac{x}{\sqrt{x^2+y^2}}u_{s\,x}(\sqrt{x^2+y^2},0,z,t)\\[10pt]
&&\dsp
u_{s\,y}(x,y,z,t)=\frac{y}{\sqrt{x^2+y^2}}u_{s\,x}(\sqrt{x^2+y^2},0,z,t)\\[10pt]
&&\dsp u_{s\,z}(x,y,z,t)=u_{s\,z}(\sqrt{x^2+y^2},0,z,t)
\end{eqnarray}
To state our results, we need the following notations and definitions:
\begin{enumerate}
\item {\bf Definition of the complex square root}. For $q_x\in\setC\backslash\setR^-$, we use the following definition of the square root $g(q_x)=q_x^{1/2}$:
\[
  g(q_x)^2=q_x \quad \hbox{ and } \quad
  \Re e[g(q_x)]>0.
\]
The branch cut of $g(q_x)$ in the complex plane will thus be the half-line defined by $\{q_x \in \setR^- \}$ (see Fig.~\ref{simple:fig:15}).
In the following, we'll use the abuse of notation $g(q_x)=\ic\sqrt{-q_x}$ for $q_x\in\setR^-$.
\\[10pt]
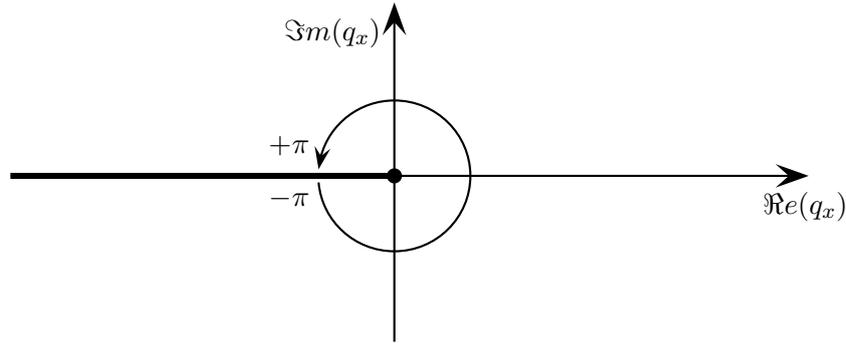
\begin{figure}[htbp]
\setlength{\unitlength}{1cm}
\psset{xunit=1cm,yunit=1cm,runit=1cm}
\begin{picture}(10,5)
\put(0,-2){\put(6.6,6.){$\Im m(q_x)$}
\put(12.9,3.7){$\Re e(q_x)$}
\put(6.4,3.8){$-\pi$}
\put(6.4,4.5){$+\pi$}
\psline[arrows=->,arrowsize=.2 4](3,4.2)(13.5,4.2)
\psline[arrows=->,arrowsize=.2 4](8.05,2)(8.05,6.5)
\psline[linewidth=0.08](3,4.2)(8.05,4.2)
\pscircle[fillstyle=solid,fillcolor=black](8.05,4.2){.1}
\psarc[arrows=->,arrowsize=.1 4](8.05,4.2){1}{-175}{175}}
\end{picture}
\caption{Definition of the function $x\mapsto (x)^{1/2}$}
\label{simple:fig:15}
\end{figure}
\item {\bf Definition of the fictitious velocities}
For a given $q\in\setR$, we define the fictitious
velocities $\mat V^+(q)$ and $\mat V_i^-(q)$ for
$i\in\{Pf,Ps,S\}$ by 
$$\mat V^+:=\mat V^+(q)=V^+\sqrt{\frac{1}{1+{V^+}^2q^2}}
\hbox{ and }
\mat V^-_i:=\mat V^-_i(q)=V^-_i\sqrt{\frac{1}{1+{V^-_i}^2q^2}}.$$
These fictitious velocities will be helpful to turn
the 3D-problem into the sum of 2D-problems indexed
by the variable $q$. Note that $\mat V^+(0)$ and
$\mat V^-_i(0)$ correspond to the real
velocities $V^+$ and $V^-_i$.
\item {\bf Definition of the functions $\kappa^{+}$
    and $\kappa^{-}_{i}$}. For $i\in\{Pf,Ps,S\}$
  and $(q_x,q_y)\in\setC\times\setR$, we define the functions
$$\kappa^+:=\kappa^+(q_x,q_y)=\left(\frac{1}{{V^+}^2}+q_x^2+q_y^2\right)^{1/2}=
\left(\frac{1}{{\mat V^+}^2(q_y)}+q_x^2\right)^{1/2}.$$
and
$$\kappa_i^-:=\kappa_i^-(q_x,q_y)=\left(\frac{1}{{V_i^-}^2}+q_x^2+q_y^2\right)^{1/2}=\left(\frac{1}{{\mat V_i^-}^2(q_y)}+q_x^2\right)^{1/2}.$$
\item {\bf Definition of the reflection and transmission coefficients.}
For a given $(q_x,q_y)\in\setC\times\setR$, we denote by $\rf(q_x,q_y)$,
 $\tf(q_x,q_y)$, $\ts(q_x,q_y)$ and 
$\tpsi(q_x,q_y)$ the solution of the linear system
\begin{equation}\label{eq:linsyst}
\mat A(q_x,q_y)
\left[
  \begin{array}{l}
    \rf(q_x,q_y) \\[10pt]
 \tf(q_x,q_y)\\[10pt]
 \ts(q_x,q_y)\\[10pt]
\tpsi(q_x,q_y)
  \end{array}
\right]
=-\frac{1}{2\ka(q_x,q_y){V^+}^2}
\left[
  \begin{array}{c}
\dsp     \frac{\ka(q_x,q_y)}{\rho^+}\\[10pt]
\dsp  1\\[10pt]
\dsp 0\\[10pt]
1
  \end{array}
\right]
,
\end{equation}
where the matrix $\mat A(q_x,q_y)$ is defined by:
$$
\begin{array}{l}
\mat A(q_x,q_y)=\\[10pt]
 \left[\begin{array}{cccc}
\dsp -\frac{\ka(q_x,q_y)}{\rho^+}&\dsp (\mat P_{11}+\mat P_{21})\kaf(q_x,q_y)&\dsp (\mat P_{12}+\mat P_{22})\kas(q_x,q_y)&\dsp \left(1-\frac{\rho_f^-}{\rho_w^-}\right)(q_x^2+q_y^2)\\[10pt]
1&\dsp \frac{m^-}{{\Vpf^-}^2}\left(\beta^-\mat P_{11}+\mat P_{21}\right)&
\dsp \frac{m^-}{{\Vps^-}^2}\left(\beta^-\mat P_{12}+\mat P_{22}\right)
&0\\[10pt]
0&2\mat P_{11}\kaf(q_x,q_y)&2\mat P_{12}\kas(q_x,q_y)& \kapsi^2(q_x,q_y)+q_x^2+q_y^2\\[10pt]
1&\mat A_{4,2}(q_x,q_y)&\mat A_{4,3}(q_x,q_y)&2\mu^-(q_x^2+q_y^2)\kapsi(q_x,q_y)
\end{array}\right],
\end{array}
$$
with 
$$A_{4,2}(q_x,q_y)=\frac{\left(\lambda^-+m^-{\beta^-}^2\right)\mat
  P_{11}
+m^-\beta^-\mat P_{21}
}{{\Vpf^-}^2}+2\mu^-{\kaf}^2(q_x,q_y)\mat P_{11},$$
$$A_{4,3}(q_x,q_y)=\frac{\left(\lambda^-+m^-{\beta^-}^2\right)\mat P_{12}+m^-\beta^-\mat P_{22}}{{\Vps^-}^2}+2\mu^-{\kas}^2(q_x,q_y)\mat P_{12}.$$
\end{enumerate}
We also denote by $V_{\max}$ the greatest velocity in
the two media: $V_{\max}=\max(V^+,\Vpf^-,\Vps^-,\Vs^-).$\\\\
We can now present the expression of the solution
to the Green Problem :
\begin{theo}
The pressure and the displacement in the top medium
are given by
$$p^+(x,0,z,t)=p_{\hbox{inc}}^+(x,z,t)+
\frac{d\xi_{\hbox{ref}}^+}{dt}(x,z,t)
\hbox{ and }
\gr
u^+(x,0,z,t)=\gr u_{\hbox{inc}}^+(x,z,\tau)\,d\tau+\gr u_{\hbox{ref}}^+(x,y,\tau)\,d\tau,$$
and the 
displacement in the bottom medium is given by
  $$\gr u_s^-(x,0,z,t)=
\gr
u^-_{Pf}(x,z,t)+\gr u^-_{Ps}(x,z,t)+\gr u^-_{S} (x,z,t)
 %
$$ where
  \begin{itemize}
  \item $p_{\hbox{inc}}^+$ and
    $u_{\hbox{inc}}^+$ are respectively
    the pressure and the displacement of the
    incident wave and satisfy~:
$$\left\{
\begin{array}{lcl}
\dsp
p_{\hbox{inc}}^+(x,z,t)&:=&\dsp\frac{\delta(t-t_0)}{4\pi{V^+}^2 r}
\\[18pt]
\dsp
u_{\hbox{inc},x}(x,z,t)&:=&\dsp\frac{xtH(t-t_0)}{4\pi{V^+}^2
  r^3\rho^+}\\[18pt]
\dsp  u_{\hbox{inc},z}(x,z,t)&:=&\dsp\frac{(z-h)tH(t-t_0)}{4\pi{V^+}^2 r^3\rho^+}\
\end{array}\right.,$$
where $\delta$ and $H$ respectively denote the usual Dirac
and Heaviside distributions. Moreover we set
$r=(x^2+(z-h)^2)^{1/2}$ and $t_0=r/V^+$ denotes the time
arrival of the incident wave at point $(x,0,z)$.
\\[5pt]
  \item $\xi_{\hbox{ref}}^+$ and
    $u_{\hbox{ref}}^+$ are respectively
    the primitive of the pressure with respect to the
    time and the displacement of the
    reflected wave and satisfy~:
$$\left\{
\begin{array}{lcl}
 \dsp \xi_{\hbox{ref}}^+(x,z,t)&=&\dsp -\int_{0}^{q_1(t)}\frac{\Im
   m\Big[\ka(\upsilon(t,q))\rf(\upsilon(t,q))\Big]}{\pi^2r \sqrt{q^2+q_0^2(t)}}\,dq,
\\[15pt]
 \dsp  u_{\hbox{ref},x}^+(x,z,t)&=&\dsp-\int_{0}^{q_1(t)} \frac{\Im
   m\Big[\ic\upsilon(t,q)\ka(\upsilon(t,q))\rf(\upsilon(t,q))\Big]}{\pi^2r\rho^+ \sqrt{q^2+q_0^2(t)}}\,dq,
\\[15pt]
\dsp u_{\hbox{ref},z}^+(x,z,t)&=&\dsp -\int_{0}^{q_1(t)}\frac{\Im
  m\Big[\ka^2(\upsilon(t,q))\rf(\upsilon(t,q))\Big]}{\pi^2r\rho^+ \sqrt{q^2+q_0^2(t)}}\,dq,
\end{array}\right.$$
$\dsp\hbox{if }t_{\hbox{h}_1} <t\leq t_{0} \hbox{ and } \frac{x}{r}>\frac{V^+}{V_{\max}}$,
$$\left\{
\begin{array}{lcl}
  \dsp  \xi_{\hbox{ref}}^+(x,z,t)&=&\dsp
  \begin{array}[t]{ll}
-\dsp\int^{q_1(t)}_{q_0(t)}\frac{\Im
   m\Big[\ka(\upsilon(t,q))\rf(\upsilon(t,q))\Big]}{\pi^2r \sqrt{q^2-q_0^2(t)}}\,dq
\\[18pt]
+\dsp \int_{0}^{q_0(t)}\frac{\Re e\Big[\ka(\gamma(t,q))\rf(\gamma(t,q))\Big]}{\pi^2r \sqrt{q_0^2(t)-q^2}}\,dq,    
  \end{array}
\\[60pt]
  \dsp u_{\hbox{ref},x}^+(x,z,t)&=&\dsp 
 \begin{array}[t]{ll}
-&\dsp \int^{q_1(t)}_{q_0(t)}\frac{\Im
   m\Big[(\ic\upsilon(t,q)\ka(\upsilon(t,q))\rf(\upsilon(t,q))\Big]}{\pi^2r\rho^+ \sqrt{q^2-q_0^2(t)}}\,dq
\\[18pt]
 +&\dsp\int_{0}^{q_0(t)}\frac{\Re
   e\Big[\ic\gamma(t,q)\ka(\gamma(t,q))\rf(\gamma(t,q))\Big]}{\pi^2r\rho^+ \sqrt{q_0^2(t)-q^2}}\,dq,
\end{array}
\\[60pt]
\dsp u_{\hbox{ref},y}^+(x,z,t)&=&\dsp 
  \begin{array}[t]{ll}
-&\dsp\int^{q_1(t)}_{q_0(t)}\frac{\Im
  m\Big[\ka^2(\upsilon(t,q))\rf(\upsilon(t,q))\Big]}{\pi^2r\rho^+ \sqrt{q^2-q_0^2(t)}}\,dq
\\[15pt]
+&\dsp \int_{0}^{q_0(t)} \frac{\Re
  e\Big[\ka^2(\gamma(t))\rf(\gamma(t))\Big]}{\pi^2r\rho^+ \sqrt{q_0^2(t)-q^2}}\,dq,
\end{array}
\end{array}\right.$$
$\dsp\hbox{if }t_{0} <t\leq t_{\hbox{h}_2}  \hbox{ and } \frac{x}{r}>\frac{V^+}{V_{\max}}$,
$$\left\{
\begin{array}{lcl}
  \dsp  \xi_{\hbox{ref}}^+(x,z,t)&=&\dsp\int^{q_0(t)}_{0}\frac{\Re e\Big[\ka(\gamma(t,q))\rf(\gamma(t,q))\Big]}{\pi^2r \sqrt{q_0^2(t)-q^2}}\,dq,
\\[15pt]
  \dsp u_{\hbox{ref},x}^+(x,z,t)&=&\dsp  \int^{q_0(t)}_{0}\frac{\Re e\Big[\ic\gamma(t,q)\ka(\gamma(t,q))\rf(\gamma(t,q))\Big]}{2\pi^2r\rho^+ \sqrt{q_0^2(t)-q^2}}\,dq,
\\[15pt]
\dsp u_{\hbox{ref},y}^+(x,z,t)&=&\dsp \int^{q_0(t)}_{0} \frac{\Re e\Big[\ka^2(\gamma(t))\rf(\gamma(t))\Big]}{\pi^2r\rho^+ \sqrt{q_0^2(t)-q^2}}\,dq,
\end{array}\right.$$
$\dsp\hbox{if } t_{\hbox{h}_2}<t
\hbox{ and } \frac{x}{r}>\frac{V^+}{V_{\max}}$ or
 $\dsp\hbox{if }t_{0} <t
\hbox{ and } \frac{x}{r}\leq \frac{V^+}{V_{\max}}$
and 
$$
\xi_{\hbox{ref}}(x,y,t)=0 \hbox{ and } \gr
u_{\hbox{ref}}(x,y,t)=0 \hbox{ else }.$$
We set here $r=(x^2+(z+h)^2)^{1/2}$  and
$t_0=r/V^+$ denotes the arrival time of the reflected
volume wave at point $(x,0,z)$, 
\begin{eqnarray}
  t_{h_1}&=&(z+h)\sqrt{\frac{1}{{V^+}^2}-\frac{1}{V^2_{\max}}}+\frac{|x|}{V_{\max}}
\end{eqnarray}
denotes the arrival time of the reflected
head-wave at point $(x,0,z)$ and 
\begin{eqnarray}
  t_{h_2}&=&\frac{r}{z+h}\sqrt{\frac{1}{{V^+}^2}-\frac{1}{{V^2_{\max}}}}
\end{eqnarray}
denotes the time after which there is no longer head
wave at point $(x,0,z)$, (contrary to the 2D case, this time does not
coincide with the arrival time of the volume wave).
We also define the functions $\gamma$, $\upsilon$,
 $q_0$ and $q_1$ by
$$\gamma : \{t\in\setR \,|\, t>t_0\}\times\setR\mapsto
\setC:=
\gamma(t,q_y)=\ic \frac
{xt}{r^2}+\frac{z+h}{r}\sqrt{\frac{t^2}{r^2}-\frac{1}{{\mat V^+}^2(q_y)}}$$
$$\upsilon : \{t\in\setR \,|\, t_{h_1}<t<t_{h_2}\}\times\setR\mapsto
\setC:=
\upsilon(t,q_y)=-\ic\left(\frac{z+h}{r}-\sqrt{\frac{1}{{\mat V^+}^2(q_y)}-\frac{t^2}{r^2}}+\frac
  {x}{r^2}t\right),$$
$$q_0:\setR\to\setR:=q_0(t)=\sqrt{\left|\frac{t^2}{r^2}-\frac{1}{
      {V^+}^2}\right|}$$
and 
$$q_1:\setR\to\setR:=q_1(t)=\sqrt{\frac{1}{x^2}\left(t-(z+h)\sqrt{\frac{1}{{V^+}^2}-\frac{1}{V^2_{\max}}}\right)^2-\frac{1}{V_{\max}^2}}.$$
\begin{remark}
  For the practical computation of the pressure, we won't have to
  explicitly compute the derivative of the function
  $\xi_{\hbox{ref}}^+$ (which would be rather tedious), since
$$p^+_{\hbox{ref}}\ast f=\partial_t\xi_{\hbox{ref}}^+\ast f=\xi_{\hbox{ref}}^+\ast f'. $$
Therefore, we'll only have to compute the derivative of the source function $f$.
\end{remark}

\item  $\gr u_{Pf}^-(x,z,t)$ is the displacement of the
      transmitted $Pf$  wave and satisfies:
$$
\left\{
\begin{array}{ll}
 \dsp u^-_{Pf,x}(x,z,t)=-\frac{\mat
   P_{11}}{\pi^2}\int_{0}^{q_{1}(t)}
\Re
 e\left[\ic\upsilon(t,q)\tf(\upsilon(t,q))\frac{\partial \upsilon}{\partial t}(t,q)\right]\,dq,
\\[18pt]
\dsp u^-_{Pf,z}(x,z,t)= 
\frac{\mat P_{11}}{\pi^2}\int_{0}^{q_{1}(t)}
\Re
e\left[\kaf(\upsilon(t,q))\tf(\upsilon(t,q))\frac{\partial\upsilon}{\partial    t}(t,q)\right]\,dq,
\end{array}\right.$$
$\dsp\hbox{if }t_{h_1} <t\leq t_{0} \hbox{ and } \left|\Im
    m\left[\gamma(t_{0},0)\right]\right|
  <\frac{1}{V_{\max}},$
$$\left\{
\begin{array}{lcl}
  \dsp u^-_{Pf,x}(x,z,t)&=&\dsp 
\begin{array}[t]{ll}
-&\dsp\frac{\mat
    P_{11}}{\pi^2}
\int_{0}^{q_0(t)}
\Re
e\left[\ic\gamma(t,q)\tf(\gamma(t,q))\frac{\partial
    \gamma}{\partial t}(t,q)\right]\,dq
\\[18pt]
-&\dsp\frac{\mat
    P_{11}}{\pi^2}
\int_{q_0(t)}^{q_{1}(t)}
\Re
 e\left[\ic\upsilon(t,q)\tf(\upsilon(t,q))\frac{\partial \upsilon}{\partial t}(t,q)\right]\,dq,
\end{array}
\\[60pt]
\dsp u^-_{Pf,z}(x,z,t)&=&\dsp 
\begin{array}[t]{ll}
&\dsp \frac{\mat
  P_{11}}{\pi^2}\int_{0}^{q_0(t)}
\Re
e\left[\kaf(\gamma(t,q))\tf(\gamma(t,q))\frac{\partial
    \gamma
}{\partial t}(t,q)\right]\,dq
\\[18pt]
+&\dsp\frac{\mat
  P_{11}}{\pi^2}
\int_{q_0(t)}^{q_{1}(t)}
\Re
e\left[\kaf(\upsilon(t,q))\tf(\upsilon(t,q))\frac{\partial\upsilon}{\partial    t}(t,q)\right]\,dq,
\end{array}
\end{array}\right.$$
$\dsp\hbox{if }t_{0} <t\leq t_{h_2}  \hbox{ and }  \left|\Im
    m\left[\gamma(t_{0},0)\right]\right|
  <\frac{1}{V_{\max}}$,
$$
\left\{
\begin{array}{ll}
  \dsp u^-_{Pf,x}(x,z,t)=- \frac{\mat P_{11}}{\pi^2}\int_{0}^{q_0(t)}
\Re
e\left[\ic\gamma(t,q)\tf(\gamma(t,q))\frac{\partial
    \gamma}{\partial t}(t,q)\right]\,dq,
\\[18pt]
\dsp u^-_{Pf,z}(x,z,t)=\frac{\mat
  P_{11}}{\pi^2}\int_{0}^{q_0(t)}
\Re
e\left[\kaf(\gamma(t,q))\tf(\gamma(t,q))\frac{\partial
    \gamma
}{\partial t}(t,q)\right]\,dq,
\end{array}\right.$$
$\dsp\hbox{if } t_{h_2}<t
\hbox{ and } \left|\Im
    m\left[\gamma(t_{0},0)\right]\right|
  <\frac{1}{V_{\max}}$ or
 $\dsp\hbox{if }t_{0} <t
\hbox{ and } \left|\Im
    m\left[\gamma(t_{0},0)\right]\right|
  \geq\frac{1}{V_{\max}}$\\
and  $\gr
 u^-_{Pf}(x,z,t)=0$ else.

$t_0$ denotes here 
 the arrival time of the $Pf$ volume wave at point $(x,0,z)$ (we recall in appendix the
  computation of $t_0$), 
\begin{eqnarray}
  t_{h_1}&=&h\sqrt{\frac{1}{{V^+}^2}-\frac{1}{V^2_{\max}}}-z\sqrt{\frac{1}{{V^-_{Pf}}^2}-\frac{1}{V^2_{\max}}}+\frac{|x|}{V_{\max}}
\end{eqnarray}
denotes the arrival time of the $Pf$
head wave at point $(x,0,z)$,
\begin{eqnarray}
  t_{h_2}&=&\frac{h^2+z^2-hz\left(\frac{c_2}{c_1}+
\frac{c_1}{c_2}
\right)+x^2}
{\frac{h}{c_1}-\frac{z}{c_2}}
\end{eqnarray}
 denotes the time after which there is no longer head
wave at point $(x,0,z)$, where
$$c_1=\sqrt{\frac{1}{{V^+}^2}-\frac{1}{{V^2_{\max}}}} \hbox{ and } c_2=\sqrt{\frac{1}{{V^2_{Pf}}}-\frac{1}{{V^2_{\max}}}}.$$
The function $q_0:[t_0\,;\,+\infty]\mapsto\setR^+$ is
the reciprocal function of $\tilde t_0
:\setR^+\mapsto:[t_0,+\infty]$, where $\widetilde
t_0(q)$ is the arrival time  at point $(x,0,z)$ of the fictitious 
$Pf$ volume wave, propagating at a velocity $\mat V^+(q)$ in the top layer and at velocity $\mat V^-_{Pf}(q)$ in the bottom layer (we recall in appendix the
  computation of $\widetilde t_0(q)$).\\[5pt]
The function $q_1:[t_1\,;\,t_0]\mapsto\setR^+$ is defined by 
$$q_{1}(t)=\sqrt{\frac{1}{x^2}\left(t+z\sqrt{\frac{1}{{V^-_{Pf}}^2}-\frac{1}{V^2_{\max}}}
-h\sqrt{\frac{1}{{V^+}^2}-\frac{1}{V^2_{\max}}}
\right)^2-\frac{1}{V_{\max}^2}}.$$
The function $\gamma:\{(t,q)\in\setR^+\times\setR^+\,|\, t>\tilde t_0(q)\}\mapsto\setC$ is implicitly defined as the only root of the function 
$${\cal
  F}(\gamma,q,t)=-z\left(\frac{1}{{\mat V_{Pf}^-}^2(q)}+\gamma^2\right)^{1/2}+h\left(\frac{1}{{\mat V^+}^2(q)}+\gamma^2\right)^{1/2}+i\gamma x-t$$
whose real part is positive.\\
The function $\upsilon:E_1\cup E_2\mapsto\setC$ is implicitly defined as the only root of the function 
$${\cal
  F}(\upsilon,q,t)=-z\left(\frac{1}{{\mat V_{Pf}^-}^2(q)}+\upsilon^2\right)^{1/2}+h\left(\frac{1}{{\mat V^+}^2(q)}+\upsilon^2\right)^{1/2}+i\upsilon x-t$$
such that $\Im m\left[\partial_t \upsilon(t,q)
\right]<0$, with
$$E_1=\left\{
(t,q)\in\setR^+\times\setR^+\,|\, t_{h_1}<t<t_{0} \hbox{ and } 0<q<q_0(t)
\right\} $$ and  $$E_2=\left\{
(t,q)\in\setR^+\times\setR^+\,|\, t_{0}<t<t_{h_1} \hbox{ and } q_0(t)<q<q_1(t)\right\}.$$
\item  $\gr u_{Ps}^-(x,z,t)$  is the displacement of the
    transmitted $Ps$ wave and satisfies:
$$
\left\{
\begin{array}{ll}
 \dsp u^-_{Ps,x}(x,z,t)=-\frac{\mat
   P_{12}}{\pi^2}\int_{0}^{q_{1}(t)}
\Re
 e\left[\ic\upsilon(t,q)\ts(\upsilon(t,q))\frac{d\upsilon}{dt}(t,q)\right]\,dq,
\\[18pt]
\dsp u^-_{Ps,z}(x,z,t)= 
\frac{\mat P_{12}}{\pi^2}\int_{0}^{q_{1}(t)}
\Re
e\left[\kas(\upsilon(t),q)\ts(\upsilon(t,q))\frac{d\upsilon}{dt}(t,q)\right]\,dq,
\end{array}\right.$$

$\dsp\hbox{if }t_{h_1} <t\leq t_{0} \hbox{ and } \left|\Im
    m\left[\gamma(t_{0},0)\right]\right|
  <\frac{1}{V_{\max}},$

$$\left\{
\begin{array}{lcl}
  \dsp u^-_{Ps,x}(x,z,t)&=&\dsp 
\begin{array}[t]{ll}
-&\dsp\frac{\mat
    P_{12}}{\pi^2}
\int_{0}^{q_0(t)}
\Re
e\left[\ic\gamma(t,q)\tf(\gamma(t,q))\frac{\partial
    \gamma}{\partial t}(t,q)\right]\,dq,
\\[18pt]
-&\dsp\frac{\mat
    P_{12}}{\pi^2}
\int_{q_0(t)}^{q_{1}(t)}
\Re
 e\left[\ic\upsilon(t,q)\ts(\upsilon(t,q))\frac{d\upsilon}{dt}(t,q)\right]\,dq
\end{array}
\\[60pt]
\dsp u^-_{Ps,z}(x,z,t)&=&\dsp 
\begin{array}[t]{ll}
&\dsp \frac{\mat
  P_{12}}{\pi^2}\int_{0}^{q_0(t)}
\Re
e\left[\kas(\gamma(t,q))\ts(\gamma(t,q))\frac{d\gamma}{dt}(t,q)\right]\,dq
\\[18pt]
+&\dsp\frac{\mat
  P_{12}}{\pi^2}
\int_{q_0(t)}^{q_{1}(t)}
\Re
e\left[\kas(\upsilon(t),q)\ts(\upsilon(t,q))\frac{d\upsilon}{dt}(t,q)\right]\,dq,
\end{array}
\end{array}\right.$$
$\dsp\hbox{if }t_{0} <t\leq t_{h_2}  \hbox{ and }  \left|\Im
    m\left[\gamma(t_{0},0)\right]\right|
  <\frac{1}{V_{\max}}$,

$$
\left\{
\begin{array}{ll}
  \dsp u^-_{Ps,x}(x,z,t)=- \frac{\mat P_{12}}{\pi^2}\int_{0}^{q_0(t)}
\Re
e\left[\ic\gamma(t,q)\ts(\gamma(t,q))\frac{d\gamma}{dt}(t,q)\right]\,dq,
\\[18pt]
\dsp u^-_{Ps,z}(x,z,t)=\frac{\mat
  P_{12}}{\pi^2}\int_{0}^{q_0(t)}
\Re
e\left[\kas(\gamma(t,q))\ts(\gamma(t,q))\frac{d\gamma}{dt}(t,q)\right]\,dq,
\end{array}\right.$$
$\dsp\hbox{if } t_{h_2}<t
\hbox{ and } \left|\Im
    m\left[\gamma(t_{0},0)\right]\right|
  <\frac{1}{V_{\max}}$ or
 $\dsp\hbox{if }t_{0} <t
\hbox{ and } \left|\Im
    m\left[\gamma(t_{0},0)\right]\right|
  \geq\frac{1}{V_{\max}}$\\
and  $\gr
 u^-_{Ps}(x,z,t)=0$ else.
$t_0$ denotes here the arrival time of the $Ps$ volume wave at point $(x,0,z)$, 
\begin{eqnarray}
  t_{h_1}&=&h\sqrt{\frac{1}{{V^+}^2}-\frac{1}{V^2_{\max}}}-z\sqrt{\frac{1}{{V^-_{Ps}}^2}-\frac{1}{V^2_{\max}}}+\frac{|x|}{V_{\max}}
\end{eqnarray}
denotes the arrival time of the $Ps$
head wave at point $(x,0,z)$ and
\begin{eqnarray}
  t_{h_2}&=&\frac{h^2+z^2-hz\left(\frac{c_2}{c_1}+
\frac{c_1}{c_2}
\right)+x^2}
{\frac{h}{c_1}-\frac{z}{c_2}} 
\end{eqnarray}
denotes the time after which there is no longer head
wave at point $(x,0,z)$, where
$$c_1=\sqrt{\frac{1}{{V^+}^2}-\frac{1}{{V^2_{\max}}}} \hbox{ and } c_2=\sqrt{\frac{1}{{V^2_{Ps}}}-\frac{1}{{V^2_{\max}}}}.$$
The function $q_0:[t_0\,;\,+\infty]\mapsto\setR^+$ is
the reciprocal function of $\tilde t_0
:\setR^+\mapsto:[t_0,+\infty]$, where $\widetilde
t_0(q)$ is the arrival time  at point $(x,0,z)$ of the fictitious 
$Ps$ volume wave, propagating at a velocity $\mat
V^+(q)$ in the top layer and at velocity $\mat 
V^-_{Ps}(q)$ in the bottom layer.\\[5pt]
The function $q_1:[t_1\,;\,t_0]\mapsto\setR^+$ is defined by 
$$q_{1}(t)=\sqrt{\frac{1}{x^2}\left(t+z\sqrt{\frac{1}{{V^-_{Ps}}^2}-\frac{1}{V^2_{\max}}}
-h\sqrt{\frac{1}{{V^+}^2}-\frac{1}{V^2_{\max}}}
\right)^2-\frac{1}{V_{\max}^2}}.$$
The function $\gamma:\{(t,q)\in\setR^+\times\setR^+\,|\, t>\tilde t_0(q)\}\mapsto\setC$ is implicitly defined as the only root of the function 
$${\cal
  F}(\gamma,q,t)=-z\left(\frac{1}{{\mat V_{Ps}^-}^2(q)}+\gamma^2\right)^{1/2}+h\left(\frac{1}{{\mat V^+}^2(q)}+\gamma^2\right)^{1/2}+i\gamma x-t$$
whose real part is positive.\\
The function $\upsilon:E_1\cup E_2\mapsto\setC$ is implicitly defined as the only root of the function 
$${\cal
  F}(\upsilon,q,t)=-z\left(\frac{1}{{\mat V_{Ps}^-}^2(q)}+\upsilon^2\right)^{1/2}+h\left(\frac{1}{{\mat V^+}^2(q)}+\upsilon^2\right)^{1/2}+i\upsilon x-t$$
such that $\Im m\left[\partial_t \upsilon(t,q) \right]<0$, with
$$E_1=\left\{
(t,q)\in\setR^+\times\setR^+\,|\, t_{h_1}<t<t_{0} \hbox{ and } 0<q<q_0(t)
\right\} $$ and  $$E_2=\left\{
(t,q)\in\setR^+\times\setR^+\,|\, t_{0}<t<t_{h_1} \hbox{ and } q_0(t)<q<q_1(t)\right\}.$$
\item  $\gr u_{S}^-(x,z,t)$  is the displacement of the
    transmitted $S$ wave and satisfies:
$$
\left\{
\begin{array}{ll}
 \dsp u^-_{S,x}(x,z,t)=-\frac{1}{\pi^2}\int_{0}^{q_{1}(t)}
\Re
 e\left[\ic \upsilon(t,q)\kapsi(\upsilon(t,q))\tpsi(\upsilon(t,q))\frac{d\upsilon}{dt}(t,q)\right]\,dq,
\\[18pt]
\dsp u^-_{S,z}(x,z,t)=
\frac{1}{\pi^2}\int_{0}^{q_{1}(t)}
\Re
e\left[(\upsilon^2(t,q)+q^2)\tpsi(\upsilon(t,q))\frac{d\upsilon}{dt}(t,q)\right]\,dq,
\end{array}\right.$$

$\dsp\hbox{if }t_{h_1} <t\leq t_{0} \hbox{ and } \left|\Im
    m\left[\gamma(t_{0},0)\right]\right|
  <\frac{1}{V_{\max}},$

$$\left\{
\begin{array}{lcl}
  \dsp u^-_{S,x}(x,z,t)&=&\dsp 
\begin{array}[t]{ll}
-&\dsp\frac{1}{\pi^2}
\int_{0}^{q_0(t)}
\Re
e\left[\ic \gamma(t,q)\kapsi(\gamma(t,q))\tpsi(\gamma(t,q))\frac{d\gamma}{dt}(t,q)\right]\,dq,
\\[18pt]
-&\dsp\frac{1}{\pi^2}
\int_{q_0(t)}^{q_{1}(t)}
\Re
 e\left[\ic \upsilon(t,q)\kapsi(\upsilon(t,q))\tpsi(\upsilon(t,q))\frac{d\upsilon}{dt}(t,q)\right]\,dq
\end{array}
\\[60pt]
\dsp u^-_{S,z}(x,z,t)&=&\dsp 
\begin{array}[t]{ll}
&\dsp \frac{1}{\pi^2}\int_{0}^{q_0(t)}
\Re
e\left[(\gamma^2(t,q)+q^2)\tpsi(\gamma(t,q))\frac{d\gamma}{dt}(t,q)\right]\,dq
\\[18pt]
+&\dsp\frac{1}{\pi^2}
\int_{q_0(t)}^{q_{1}(t)}
\Re
e\left[(\upsilon^2(t,q)+q^2)\tpsi(\upsilon(t,q))\frac{d\upsilon}{dt}(t,q)\right]\,dq,
\end{array}
\end{array}\right.$$
$\dsp\hbox{if }t_{0} <t\leq t_{h_2}  \hbox{ and }  \left|\Im
    m\left[\gamma(t_{0},0)\right]\right|
  <\frac{1}{V_{\max}}$,

$$
\left\{
\begin{array}{ll}
  \dsp u^-_{S,x}(x,z,t)= -\frac{1}{\pi^2}\int_{0}^{q_0(t)}
\Re
e\left[\ic \gamma(t,q)\kapsi(\gamma(t,q))\tpsi(\gamma(t,q))\frac{d\gamma}{dt}(t,q)\right]\,dq,
\\[18pt]
\dsp u^-_{S,z}(x,z,t)=\frac{1}{\pi^2}\int_{0}^{q_0(t)}
\Re
e\left[(\gamma^2(t,q)+q^2)\tpsi(\gamma(t,q))\frac{d\gamma}{dt}(t,q)\right]\,dq,
\end{array}\right.$$
$\dsp\hbox{if } t_{h_2}<t
\hbox{ and } \left|\Im
    m\left[\gamma(t_{0},0)\right]\right|
  <\frac{1}{V_{\max}}$ or
 $\dsp\hbox{if }t_{0} <t
\hbox{ and } \left|\Im
    m\left[\gamma(t_{0},0)\right]\right|
  \geq\frac{1}{V_{\max}}$
and\\  $\gr
 u^-_{Ps}(x,z,t)=0$ else.
$t_0$ denotes here the arrival time of the $S$ volume wave at point $(x,0,z)$ (we recall in appendix the
  computation of $t_0$), 
\begin{eqnarray}
  t_{h_1}&=&h\sqrt{\frac{1}{{V^+}^2}-\frac{1}{V^2_{\max}}}-z\sqrt{\frac{1}{{V^-_{S}}^2}-\frac{1}{V^2_{\max}}}+\frac{|x|}{V_{\max}}
\end{eqnarray}
 denotes the arrival time of the $S$
head-wave at point $(x,0,z)$ and 
\begin{eqnarray}
  t_{h_2}&=&\frac{h^2+z^2-hz\left(\frac{c_2}{c_1}+
\frac{c_1}{c_2}
\right)+x^2}
{\frac{h}{c_1}-\frac{z}{c_2}}
\end{eqnarray}
denotes the time after which there is no longer head
wave at point $(x,0,z)$, where
$$c_1=\sqrt{\frac{1}{{V^+}^2}-\frac{1}{{V^2_{\max}}}} \hbox{ and } c_2=\sqrt{\frac{1}{{V^2_{S}}}-\frac{1}{{V^2_{\max}}}}.$$
The function $q_0:[t_0\,;\,+\infty]\mapsto\setR^+$ is
the reciprocal function of $\tilde t_0
:\setR^+\mapsto:[t_0,+\infty]$, where $\widetilde
t_0(q)$ is the arrival time  at point $(x,0,z)$ of the fictitious
$S$ volume wave, propagating at a velocity $
\mat V^+(q)$ in the top layer and at velocity $\mat 
V^-_{S}(q)$ in the bottom layer (we recall in appendix the
  computation of $\widetilde t_0(q)$).\\[5pt]
The function $q_1:[t_1\,;\,t_0]\mapsto\setR^+$ is defined by 
$$q_{1}(t)=\sqrt{\frac{1}{x^2}\left(t+z\sqrt{\frac{1}{{V^-_{S}}^2}-\frac{1}{V^2_{\max}}}
-h\sqrt{\frac{1}{{V^+}^2}-\frac{1}{V^2_{\max}}}
\right)^2-\frac{1}{V_{\max}^2}}.$$
The function $\gamma:\{(t,q)\in\setR^+\times\setR^+\,|\, t>\tilde t_0(q)\}\mapsto\setC$ is implicitly defined as the only root of the function 
$${\cal
  F}(\gamma,q,t)=-z\left(\frac{1}{{\mat V_{S}^-}^2(q)}+\gamma^2\right)^{1/2}+h\left(\frac{1}{{\mat V^+}^2(q)}+\gamma^2\right)^{1/2}+i\gamma x-t$$
whose real part is positive.\\
The function $\upsilon:E_1\cup E_2\mapsto\setC$ is implicitly defined as the only root of the function 
$${\cal
  F}(\upsilon,q,t)=-z\left(\frac{1}{{\mat V_{S}^-}^2(q)}+\upsilon^2\right)^{1/2}+h\left(\frac{1}{{\mat V^+}^2(q)}+\upsilon^2\right)^{1/2}+i\upsilon x-t$$
such that $\Im m\left[\partial_t \upsilon(t,q) \right]<0$, with
$$E_1=\left\{
(t,q)\in\setR^+\times\setR^+\,|\, t_{h_1}<t<t_{0} \hbox{ and } 0<q<q_0(t)
\right\} $$ and  $$E_2=\left\{
(t,q)\in\setR^+\times\setR^+\,|\, t_{0}<t<t_{h_1} \hbox{ and } q_0(t)<q<q_1(t)\right\}.$$
  \end{itemize}
\end{theo}

\section{Proof of the theorem}
\label{sec:proof-theorem}
To prove the theorem, we use the Cagniard-de Hoop method~(see~\cite{Cag, DH,
  VDH, PG,QG}), which consists of three steps:
\begin{enumerate}
\item We apply a Laplace transform  in time, 
$$
\tilde u(x,y,z,s)=\int_0^{+\infty} u(x,y,z,t)\, e^{-st}\,dt,
$$
and a Fourier
transform in the $x$ and $y$ variables, 
$$
\hat u(k_x,k_y,z,s)=\int_{-\infty}^{+\infty}\int_{-\infty}^{+\infty}\tilde u(x,y,s)\,e^{\ic (k_xx+k_yy)}\,dx\,dy
$$ to~\eqref{equ:diagsyst} in order to obtain an ordinary differential system whose solution $\hat {\mat G}(k_x,k_y,z,s)$ can be explicitly computed (\S~\ref{sec:solut-lapl-four});
\item we apply an inverse Fourier transform in the
  $x$ and $y$ variables to $\mat G$ (we recall that
  we only need the solution at $y=0$:
 $$
\tilde
{\mat G}(x,0,z,s)=\frac1{4\pi^2}\int_{-\infty}^{+\infty}\int_{-\infty}^{+\infty}\hat
{\mat G}(k_x,k_y,z,s)\,e^{-\ic k_xx}\,dk_x\,dk_y.
$$
And, using tools of complex analysis, we turn the
inverse Fourier transform in the $x$ variable 
into the Laplace transform of
some function $\mat H(x,k_y,z,t)$ (\S~\ref{sec:lapl-transf-solut}):
\begin{equation}
  \label{acousporo3d2:eq:1}
\tilde {\mat G}(x,0,z,s)=\frac1{4\pi^2}\int_{-\infty}^{+\infty}\int_{0}^{+\infty}\mat H(x,k_y,z,t)\,e^{-st}\,dt\,dk_y;
  \end{equation}
\item the last step of the method consists in
  inverting the order of integration
  in~(\ref{acousporo3d2:eq:1}) to obtain
$$\tilde
{\mat G}(x,0,z,s)=\frac1{4\pi^2}\int_{0}^{+\infty}
\left(\int_{-q(t)}^{+q(t)}\mat H(x,k_y,z,t)\,dk_y\right)\,e^{-st}\,dt.$$
 Then, using the
injectivity of the Laplace transform, we identify
$\mat G(x,0,z,t)$ to 
$$\frac1{4\pi^2}\int_{-q(t)}^{+q(t)}\mat H(x,k_y,z,t)\,dk_y$$
 (see \S~\ref{sec:inv_integ}).
\end{enumerate}
\subsection{The solution in the Laplace-Fourier plane}
\label{sec:solut-lapl-four}
Let us  first
apply a Laplace transform in time 
and a Fourier
transform in the $x$ and $y$ variables
to~\eqref{equ:diagsyst} to obtain 
\begin{equation}\label{de-rr:eq:3}
\left\{
\begin{array}{lll}
\dsp \left(\frac{s^2}{{V^+}^2}+k_x^2+k_y^2\right )\hat p^+-\frac{
  \partial^2\hat p^+}{\partial z^2}=\frac{\delta(z-h)}{{V^+}^2},&&y>0,\\[18pt]
\dsp \left(\frac{s^2}{{V_i^-}^2}+k_x^2+k_y^2\right)\hat\Phi_{i}^--\frac{ \partial^2\hat\Phi_{i}^-}{\partial z^2}=0,\;\;i\in\{Pf,Ps,S\}&&y<0,\\[18pt]
{\hat{\cal B}} (\hat p^+,\hat\Phi^-_{Pf},\hat \Phi^-_{Pf},\hat\Phi_S^-)=0&& y=0,
\end{array}
\right.
\end{equation}
where $\hat{\cal B}$ is the Laplace-Fourier
transform of the operator $\cal B$.\\\\ 
From the two first equations
of~(\ref{de-rr:eq:3}), we deduce that the solution
$(\hat p^+, (\hat\Phi_{i}^-)_{i\in\{Pf,Ps,S\}})$
is such that
\begin{equation}
\left\{
  \begin{array}{lll}
\label{acousporo:eq:5}
 \dsp \hat p^+= \hat p_{\hbox{inc}}^+ +\hat p_{\hbox{ref}}^+,\\[10pt]
\dsp \hat p_{\hbox{inc}}^+=\frac{1}{s{V^+}^2\ka\left(\frac{k_x}s,\frac{k_y}s\right)}e^{-s|z-h|\ka\left(\frac{k_x}s,\frac{k_y}s\right)},
\dsp \hat p_{\hbox{ref}}^+=R(k_x,k_y,s)e^{-sz\ka\left(\frac{k_x}s,\frac{k_y}s\right)},\\[20pt]
\dsp \hat \Phi^-_{i}=T_{i}(k_x,k_y,s)e^{-s\left(z\kappa_i^-\left(\frac{k_x}s\frac{k_y}s\right)\right)}, \quad i\in\{Pf,Ps,S\},\\[10pt]
  \end{array}
\right.
\end{equation}
where the coefficients $R$ and $T_{i}$ are
computed by using the last equation
of~(\ref{de-rr:eq:3}):
$$
{\hat{\cal B}} ( \hat p_{\hbox{ref}}^+,\hat\Phi^-_{Pf},\hat \Phi^-_{Ps},\hat\Phi^-_S)=-{\hat{\cal B}} (\hat p_{\hbox{inc}}^+,0,0,0),
$$
or, from(~\ref{acousporo:eq:5}):
$$
{\hat{\cal B}}
\left[R(k_x,k_y,s),T_{Pf}(k_x,k_y,s),T_{Ps}(k_x,k_y,s),T_S(k_x,k_y,s)\right]=-{\hat{\cal
    B}}
\left(\frac{e^{-sh\ka\left(\frac{k_x}s,\frac{k_y}s\right)}}{s\ka\left(\frac{k_x}s,\frac{k_y}s\right)},0,0,0\right).
$$
After some calculations that we don't detail here,
we obtain that $R(k_x,k_y,s)$,
$T_{Pf}(k_x,k_y,s)$, $T_{Ps}(k_x,k_y,s)$, and $T_S(k_x,k_y,s))$ are
solution to
\begin{equation}\label{acousporo:eq:9}
\mat A\left(\frac{k_x}s,\frac{k_y}s\right)
\left[
  \begin{array}{l}
   R(k_x,k_y,s) \\[10pt]
s^2 T_{Pf}(k_x,k_y,s)\\[10pt]
s^2 T_{Ps}(k_x,k_y,s)\\[10pt]
s^3 T_S(k_x,k_y,s)
  \end{array}
\right]
=-\frac{e^{-sh\ka\left(\frac{k_x}s,\frac{k_y}s\right)}}{2s\ka\left(\frac{k_x}s,\frac{k_y}s\right){V^+}^2}
\left[
  \begin{array}{c}
\dsp     \frac{\ka\left(\frac{k_x}s,\frac{k_y}s\right)}{\rho^+}\\[10pt]
\dsp  1\\[10pt]
\dsp 0\\[10pt]
1
  \end{array}
\right]
.
\end{equation}
From the definition of the reflection and
transmission coefficients we deduce that
\begin{equation}\label{acousporo:eq:10}
\left[
  \begin{array}{l}
   R(k_x,k_y,s) \\[10pt]
s^2 T_{Pf}(k_x,k_y,s)\\[10pt]
s^2 T_{Ps}(k_x,k_y,s)\\[10pt]
s^3 T_S(k_x,k_y,s)
  \end{array}
\right]
=\frac 1 s \left[
  \begin{array}{l}
    \rf\left(\frac{k_x}s,\frac{k_y}s\right) \\[10pt]
 \tf\left(\frac{k_x}s,\frac{k_y}s\right)\\[10pt]
 \ts\left(\frac{k_x}s,\frac{k_y}s\right)\\[10pt]
\tpsi\left(\frac{k_x}s,\frac{k_y}s\right)
  \end{array}
\right]e^{-sh\ka\left(\frac{k_x}s\right)}
.
\end{equation}

Finally, we obtain:
\begin{equation}
\left\{
  \begin{array}{lll}
    \label{de-rr:eq:12}
 \dsp \hat p^+= \hat p_{\hbox{inc}}^+ +\hat p_{\hbox{ref}}^+,\\[10pt]
\dsp \hat p_{\hbox{inc}}^+=\frac{1}{s{V^+}^2\ka\left(\frac{k_x}s,\frac{k_y}s\right)}e^{-s|y-h|\ka\left(\frac{k_x}s,\frac{k_y}s\right)},
\dsp \hat p_{\hbox{ref}}^+=\frac 1 s \rf\left(\frac{k_x}s,\frac{k_y}s\right)e^{-s(z+h)\ka\left(\frac{k_x}s,\frac{k_y}s\right)},\\[22pt]
\dsp \hat \Phi^-_{i}=\frac{1}{s^3}{\cal
  T}_{i}\left(\frac{k_x}{s},\frac{k_y}s\right)e^{-s\left(z\kappa_i^-\left(\frac{k_x}s,\frac{k_y}s\right)-h\kappa^+\left(\frac{k_x}s,\frac{k_y}s\right)\right)}, \quad i\in\{Pf,Ps\},\\[18pt]
\dsp \hat \Phi^-_{S}=\frac{1}{s^4}{\cal T}_{S}\left(\frac{k_x}{s},\frac{k_y}s\right)e^{-s\left(z\kappa_S^-\left(\frac{k_x}s,\frac{k_y}s\right)-h\kappa^+\left(\frac{k_x}s,\frac{k_y}s\right)\right)}.
  \end{array}
\right.
\end{equation}
and 
\begin{equation}
\left\{
  \begin{array}{lll}
\label{acousporo:eq:1}
\dsp \hat u^+= \hat u_{\hbox{inc}}^+ +\hat u_{\hbox{ref}}^+,\\[10pt]
\dsp \hat u^+_{\hbox{inc},x}=\ic\frac{k_x}{\ro^+s^2} \hat p_{\hbox{inc}}^+,\quad \dsp \hat u^+_{\hbox{inc},z}=\hbox {sign}(h-z)\frac{{\ka}\left(\frac{k_x}s,\frac{k_y}s\right)}{\ro^+s}\hat p_{\hbox{inc}}^+,\\[10pt]
\dsp \hat u^+_{\hbox{ref},x}=\ic\frac{k_x}{\ro^+s^2} \hat p_{\hbox{ref}}^+,\quad \dsp \hat u^+_{\hbox{ref},z}=\frac{{\ka}\left(\frac{k_x}s,\frac{k_y}s\right)}{\ro^+s}\hat p_{\hbox{ref}}^+,\\[12pt]
\dsp \hat u^-_{sx}=-\ic k_x \mat P_{11}\hat \Phi_{Pf}-\ic k_x\mat P_{12}  \hat \Phi_{Ps}^--isk_x\kapsi\left(\frac{k_x}s\right)\Phi_S^-\\[20pt]
\dsp \hat u^-_{sz}=s\kaf\left(\frac{k_x}s,\frac{k_y}s\right)\mat P_{11} \hat \Phi_{Pf}^-+s\kas\left(\frac{k_x}s,\frac{k_y}s\right)\mat P_{12} \hat \Phi_{Ps}^-+(k_x^2+k_y^2)\Phi_S^-
  \end{array}
\right.
\end{equation}
In the following we only
detail the computation of $\hat u^-_{sx,Ps}=-\ic
k_x \mat P_{12}\,\hat \Phi_{Ps}^-$, since the
computation of the other terms
is  very similar. 
\subsection{The Laplace transform of the solution}
\label{sec:lapl-transf-solut}
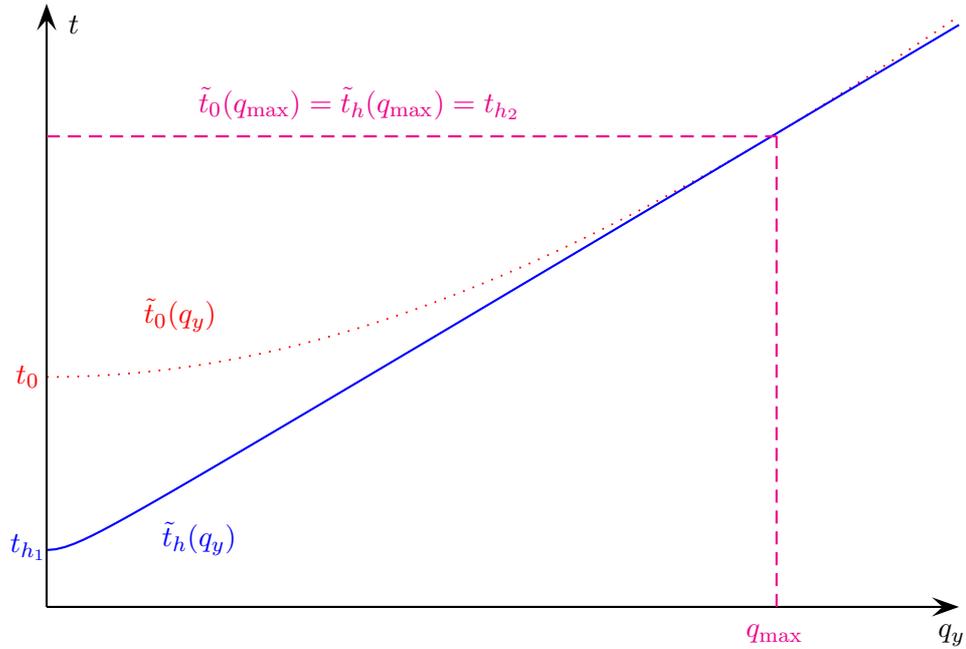
\begin{figure}[htbp]
\setlength{\unitlength}{.8cm}
\psset{xunit=0.8cm,yunit=0.8cm}
\begin{picture}(15,11)(-3,-.5)
\put(1,0){\psplot[linecolor=red,plotpoints=100,linestyle=dotted]{0}{15}{x 2 mul 2 exp 1000 add sqrt .5 mul -12 add}
\psplot[linecolor=blue,plotpoints=100]{0}{15}{x 2 mul 2 exp 1 add sqrt .6 mul 999 sqrt 0.8 mul add .5 mul -12 add}
\psline[arrowsize=.2 2]{->}(0,0)(0,10)
\psline[arrowsize=.2 2]{->}(0,0)(15,0)
\put(0.2,9.5){ $t$}
\put(14.5,-0.5){ $q_y$}
\put(11.5,-0.5){\magenta $q_{\max}$}
\put(-.5,3.7){\red $t_{0}$}
\put(-.6,.9){\blue $t_{h_1}$}
\put(1.6,4.7){\red $\tilde t_{0}(q_y)$}
\put(1.9,1.){\blue $\tilde t_{h}(q_y)$}
\put(2.5,8.2){\magenta  $\tilde t_{0}(q_{\max})=\tilde t_{h}(q_{\max})=t_{h_2}$}
\psline[linestyle=dashed,linecolor=magenta](12,0)(12,7.8)
\psline[linestyle=dashed,linecolor=magenta](0,7.8)(12,7.8)}
\end{picture}
  \caption{\label{acousporo3d:fig:1} Functions
    $\tilde t_0$  (red,  dotted line) and
    $\tilde t_{h}$ (blue, solid line)}
\end{figure}

\noindent  We apply an inverse Fourier transform in the
$x$ and $y$ variable to $\hat u^-_{sx,Ps}$ and we
set $k_x=q_xs$  and  $k_y=q_ys$ to obtain (we
recall that we consider $y=0$) 
  \begin{eqnarray*}
 \tilde
 u^-_{sx,Ps}(x,0,z,s)\!\!\!\!&=&\!\!\!\!-\int_{-\infty}^{+\infty}\int_{-\infty}^{+\infty}\frac{\ic q_x\mat P_{12}}{4\pi^2}\ts(q_x,q_y)e^{-s\left(-z\kas(q_x,q_y)+h\ka(q_x,q_y))+\ic q_xx\right)}\,dq_x\,dq_y\\[10pt]
&=&\!\!\!\!-\frac{\mat P_{12}}{4\pi^2}\int_{-\infty}^{+\infty}\int_{-\infty}^{+\infty}\Xi(q_x,q_y)\,dq_x\,dq_y,
  \end{eqnarray*}
with 
$$
\Xi(q_x,q_y)=\ic q_x\,\ts(q_x,q_y)e^{-s\left(-z\kas(q_x,q_y)+h\ka(q_x,q_y))+\ic q_xx\right)}.
$$
Let us now focus on the integral over $q_x$ for a
fixed $q_y$
  \begin{eqnarray}
\label{acousporo3d:eq:1}
\int_{-\infty}^{+\infty}\Xi(q_x,q_y)\,dq_x&=&\int_{-\infty}^{+\infty}\ic q_x\ts(q_x,q_y)e^{-s\left(-z\kas(q_x,q_y)+h\ka(q_x,q_y))+\ic q_xx\right)}\,dq_x
  \end{eqnarray}
This integral is very similar to the one we have
obtained in 2D~\cite{RAP_DE6509}, therefore, using the
same method, we have:
\begin{itemize}
 \item if $\dsp
   |\gamma(q_y,t_{0})|>\frac{1}{\mat V_{\max}(q_y)}$
\begin{equation*}
  \int_{-\infty}^{+\infty}\Xi(q_x,q_y)\,dq_x=2
  \int_{\tilde t_{0}(q_y)}^{+\infty} 
\Re e
\left(\ic\gamma(q_y,t)\ts(q_y,\gamma(q_y,t))\frac{\partial
    \gamma(q_y,t)}{\partial t}\right)
 e^{-st}dt
\end{equation*}
 \item if $\dsp
   |\gamma(q_y,\tilde t_0(q_y))|<\frac{1}{\mat V_{\max}(q_y)}$
$$\begin{array}{rcl}
\dsp \int_{-\infty}^{+\infty}\Xi(q_x,q_y)\,dq_x&=&\dsp 2\int_{\tilde t_h(q_y)}^{\tilde t_{0}(q_y)}{\Re e\left(
 \ic\upsilon(q_y,t)\ts(q_y,\upsilon(q_y,t))\frac{\partial\upsilon(q_y,t)}{\partial t}\right)e^{-st}dt
}\\[16pt]
&+&2\dsp\int_{\tilde t_0(q_y)}^{+\infty}{\Re e\left(
 \ic\gamma(q_y,t)\ts(q_y,\gamma(q_y,t))\frac{\partial \gamma(q_y,t)}{\partial t}\right)e^{-st}dt
}
\end{array}
$$
\end{itemize}
where $\mat V_{\max}$ is greatest fictitious velocity defined by:
$$\mat V_{\max}:=\mat V_{\max}(q)=V_{\max}\sqrt{\frac{1}{1+{V_{\max}}^2q^2}},$$
$\tilde t_0$ is the fictitious arrival time of the
$Ps$ volume wave we have defined in the theorem,
$\tilde t_h(q_y)$ the fictitious arrival time of the $Ps$
head wave defined by
$$\tilde t_h:= \tilde t_h(q)=h\sqrt{\frac{1}{{\mat V^+}^2(q)}-\frac{1}{\mat V_{\max}^2(q)}}-z\sqrt{\frac{1}{{\mat V_{Ps}^-}^2(q)}-\frac{1}{\mat V_{\max}^2(q)}}+\frac{|x|}{\mat V_{\max}(q)}.$$
Let us recall~\cite{diaz_th} that the condition $\dsp
   |\gamma(q_y,\tilde t_0(q_y))|<\frac{1}{\mat V_{\max}(q_y)}$ is equivalent to 
$$|\gamma(0,t_0)|<\frac 1{V_{\max}}
\hbox{ and } |q_y|\leq q_{\max},$$
with 
$$q_{\max}=\sqrt{\frac{r^2}{\left(
\frac{h}{\sqrt{\frac{1}{{V^+}^2}-\frac{1}{{V_{\max}}^2} }}
-\frac{z}{\sqrt{\frac{1}{{V^-_{Ps}}^2}-\frac{1}{{V_{\max}}^2} }}
\right)}-\frac{1}{{V_{\max}}^2}}.$$
Moreover,
  $\tilde t_{h}$ is bijective from
  $[0\,;\,q_{\max}]$ to
  $[t_0\,;\,\tilde t_{h}(q_{\max})]$
  and we denote its inverse by $q_{h}$:
$$q_{h}(t)=\sqrt{\frac{1}{x^2}\left(t+z\sqrt{\frac{1}{{V^-_{Ps}}^2}-\frac{1}{V^2_{\max}}}
-h\sqrt{\frac{1}{{V^+}^2}-\frac{1}{V^2_{\max}}}
\right)^2-\frac{1}{V_{\max}^2}}.$$
Let us also recall that for $q_y=q_{\max}$, the arrival times
   of the fictitious head and volume waves are the same: 
$\tilde t_{h}(q_{\max})=\tilde t_{0}(q_{\max}).$
As an illustration, we represent the functions
$\tilde t_{0}$ (the red dotted line) and
$\tilde t_{h}$ (the blue solid line) in Fig.~\ref{acousporo3d:fig:1}.

We then deduce that
\begin{itemize}
\item if $|\gamma(0,t_{0})|\geq \frac
  1{V_{\max}}$
\begin{equation*}
  \int_{\RR^2}\Xi(q_x,q_y)\,dq_xdq_y=
  2\int_{\RR} \int_{\tilde t_{0}(q_y)}^{+\infty} 
\Re e
\left(\ic\gamma(q_y,t)\ts(q_y,\gamma(q_y,t))\frac{\partial
    \gamma(q_y,t)}{\partial t}\right)
 e^{-st}\,dtdq_y;
\end{equation*}
\item if $|\gamma(0,t_{0})|< \frac
  1{V_{\max}}$
$$\hspace*{-0.7cm}\begin{array}{rcl}
\dsp \int_{\RR^2}\Xi(q_x,q_y)\,dq_xdq_y&=& 2\dsp\int_{-q_{\max}}^{+q_{\max}} \int_{\tilde t_{h}(q_y)}^{\tilde t_{0}(q_y)}{\Re e\left(
 \ic\upsilon(q_y,t)\ts(q_y,\upsilon(q_y,t))\frac{\partial\upsilon(q_y,t)}{\partial t}\right)e^{-st}\,dtdq_y
}\\[16pt]
&+&\dsp 2\int_{\RR}\int_{\tilde t_{0}(q_y)}^{+\infty}{\Re e\left(
 \ic\gamma(q_y,t)\ts(q_y,\gamma(q_y,t))\frac{\partial \gamma(q_y,t)}{\partial t}\right)e^{-st}\,dtdq_y.
}
\end{array}
$$
\end{itemize}
\subsection{Inversion of the integrals}
\label{sec:inv_integ}
The key point of the method is the inversion of
the integral with respect to $q_y$ with the
integral with respect to $t$. For the volume wave we have (see Figs.~\ref{acousporo3d:fig:3} and~\ref{acousporo3d:fig:2}), after having remark that the integrand is even with respect to $q_y$:
$$
\begin{array}{ll}
 &\dsp 
  \int_{-\infty}^{+\infty} \int_{\tilde t_{0}(q_y)}^{+\infty} 
\Re e
\left(\ic\gamma(q_y,t)\ts(q_y,\gamma(q_y,t))\frac{\partial
    \gamma(q_y,t)}{\partial t}\right)
 e^{-st}\,dtdq_y\\[16pt]
=&
\dsp   2\int_{t_{0}}^{+\infty} \int_{0}^{q_0(t)} 
\Re e
\left(\ic\gamma(q_y,t)\ts(q_y,\gamma(q_y,t))\frac{\partial
    \gamma(q_y,t)}{\partial t}\right)
 e^{-st}\,dq_ydt;
\end{array}$$
 and for the head wave (see Figs.~\ref{acoustique23d:fig:2} and~\ref{acoustique23d:fig:5}):
$$
\begin{array}{ll}
 \dsp 
 &\dsp\int_{-q_{\max}}^{+q_{\max}} \int_{\tilde t_{h}(q_y)}^{\tilde t_{0}(q_y)}\Re e\left(
 \ic\upsilon(q_y,t)\ts(q_y,\upsilon(q_y,t))\frac{\partial\upsilon(q_y,t)}{\partial t}\right)e^{-st}\,dtdq_y
\\[16pt]
=
&\dsp2\int_{t_{h_1}}^{t_0} \int_{0}^{q_h(t)}\Re e\left(
 \ic\upsilon(q_y,t)\ts(q_y,\upsilon(q_y,t))\frac{\partial\upsilon(q_y,t)}{\partial t}\right)e^{-st}\,dq_ydt\\[16pt]
+
&\dsp2\int_{t_0}^{t_{h_2}} \int_{q_0(t)}^{q_h(t)}\Re e\left(
 \ic\upsilon(q_y,t)\ts(q_y,\upsilon(q_y,t))\frac{\partial\upsilon(q_y,t)}{\partial t}\right)e^{-st}\,dq_ydt.
\end{array}$$

\setlength{\unitlength}{1cm}
\begin{figure}[H]
  \begin{tabular}{cc}
\begin{minipage}{.48\linewidth}
\setlength{\unitlength}{.4cm}
\psset{xunit=0.4cm,yunit=0.4cm}
\begin{picture}(15,11)(-1,0)
\put(1,0){\psplot[linecolor=red,plotpoints=100]{0}{15}{x 2 mul 2 exp 1000 add sqrt .5 mul -12 add}
\psline[arrowsize=.2 2]{->}(0,0)(0,10)
\psline[arrowsize=.2 2]{->}(0,0)(15,0)
\put(0.2,9.5){ $t$}
\put(14.5,-0.8){ $q_y$}
\put(-1,3.7){\red  $t_0$}
\put(1.6,2.7){\red  $\tilde t_0(q_y)$}
\psline(1,9)(1,3.8)
\psline[arrowsize=.1 1]{->}(1,3.9)(1,6.9)
\psline(2,9)(2,3.9)
\psline[arrowsize=.1 1]{->}(2,4)(2,6.9)
\psline(3,9)(3,4)
\psline[arrowsize=.1 1]{->}(3,4.3)(3,7)
\psline(4,9)(4,4.3)
\psline[arrowsize=.1 1]{->}(4,4.3)(4,7.2)
\psline(4,9)(4,4.3)
\psline(5,9)(5,4.6)
\psline[arrowsize=.1 1]{->}(5,4.9)(5,7.4)
\psline(6,9)(6,4.9)
\psline[arrowsize=.1 1]{->}(6,5.3)(6,7.6)
\psline(7,9)(7,5.3)
\psline[arrowsize=.1 1]{->}(7,5.7)(7,7.9)
\psline(8,9)(8,5.7)
\psline[arrowsize=.1 1]{->}(8,5.7)(8,8.2)
\psline(9,9)(9,7)
\psline[arrowsize=.1 1]{->}(9,6.2)(9,8.4)
\psline(10,9)(10,6.8)
\psline[arrowsize=.1 1]{->}(10,7)(10,8.7)
\psline(11,9)(11,7.2)
\psline[arrowsize=.1 1]{->}(11,7.2)(11,8.9)
\psline(12,9)(12,7.8)
\psline[arrowsize=.1 1]{->}(12,7.8)(12,9)
\psline(13,9)(13,8.5)
\psline[arrowsize=.1 1]{->}(0,4)(0,6.9)}
\end{picture}
  \caption{\label{acousporo3d:fig:3} Integration
    first over $q_y$ then over $t$ for the volume
    wave}
\end{minipage}
&\begin{minipage}{.48\linewidth}
\setlength{\unitlength}{.4cm}
\psset{xunit=0.4cm,yunit=0.4cm}
\begin{picture}(15,11)(-1,0)
\put(1,0){\psplot[linecolor=red,plotpoints=100]{0}{15}{x 2 mul 2 exp 1000 add sqrt .5 mul -12 add}
\psline[arrowsize=.2 2]{->}(0,0)(0,10)
\psline[arrowsize=.2 2]{->}(0,0)(15,0)
\put(0.2,9.5){ $t$}
\put(14.5,-0.8){ $q_y$}
\put(-1,3.7){\red   $t_0$}
\put(3.6,3.5){\red  $q_0(t)$}
 \psline(0,5)(6.2,5)
 \psline[arrowsize=.1 1]{->}(0,5)(2.8,5)
 \psline(0,6)(8.4,6)
 \psline[arrowsize=.1 1]{->}(0,6)(4,6)
 \psline(0,7)(10.4,7)
 \psline[arrowsize=.1 1]{->}(0,7)(5.5,7)
 \psline(0,8)(12.3,8)
 \psline[arrowsize=.1 1]{->}(0,8)(6.5,8)
 \psline(0,9)(13.8,9)
 \psline[arrowsize=.1 1]{->}(0,9)(7.5,9)
}
\end{picture}
  \caption{\label{acousporo3d:fig:2}Integration
    first over $t$ then over $q_y$ for the volume wave}
\end{minipage}
\end{tabular}
\end{figure}

\setlength{\unitlength}{1cm}
\begin{figure}[H]
  \begin{tabular}{cc}
\begin{minipage}{.48\linewidth}
\setlength{\unitlength}{.4cm}
\psset{xunit=0.4cm,yunit=0.4cm}
\begin{picture}(15,11)(-1,0)
\put(1,0){\psplot[linecolor=red,plotpoints=100]{0}{15}{x 2 mul 2 exp 1000 add sqrt .5 mul -12 add}
\psplot[linecolor=blue,plotpoints=100]{0}{15}{x 2 mul 2 exp 1 add sqrt .6 mul 999 sqrt 0.8 mul add .5 mul -12 add}
\psline[arrowsize=.2 2]{->}(0,0)(0,10)
\psline[arrowsize=.2 2]{->}(0,0)(15,0)
\put(0.5,9.5){$ t$}
\put(14.5,-0.8){$q_y$}
\put(11.,-0.7){\magenta  $q_{\max}$}
\put(-1,3.7){\red $t_{0}$}
\put(-1.2,.9){\blue$t_{h_1}$}
\put(1.6,4.7){\red  $\tilde t_{0}(q_y)$}
\put(1.9,1.){\blue  $\tilde t_{h}(q_y)$}
\psline(1,1.3)(1,3.8)
\psline[arrowsize=.1 1]{->}(1,1.3)(1,2.9)
\psline(2,1.9)(2,3.9)
\psline[arrowsize=.1 1]{->}(2,2.3)(2,3.1)
\psline(3,2.5)(3,4)
\psline[arrowsize=.1 1]{->}(3,2.8)(3,3.5)
\psline(4,3.1)(4,4.3)
\psline[arrowsize=.1 1]{->}(4,3.2)(4,3.9)
\psline(4,3.1)(4,4.3)
\psline[arrowsize=.1 1]{->}(4,3.2)(4,3.9)
\psline(5,3.7)(5,4.6)
\psline[arrowsize=.1 1]{->}(5,4.2)(5,4.3)
\psline(6,4.25)(6,4.9)
\psline[arrowsize=.1 1]{->}(6,4.5)(6,4.75)
\psline(7,4.9)(7,5.3)
\psline[arrowsize=.1 1]{->}(7,4.99)(7,5.3)
\psline(8,5.45)(8,5.7)
\psline[linestyle=dashed,linecolor=magenta](12,0)(12,7.8)
\psline[arrowsize=.1 1]{->}(0,1)(0,2.8)}
\end{picture}
  \caption{\label{acoustique23d:fig:2} Integration
    first over $q_y$ then over $t$ for the head
    wave}
\end{minipage}
&\begin{minipage}{.48\linewidth}
\setlength{\unitlength}{.4cm}
\psset{xunit=0.4cm,yunit=0.4cm}
\begin{picture}(15,11)(-1,0)
\put(1,0){\psplot[linecolor=red,plotpoints=100]{0}{15}{x 2 mul 2 exp 1000 add sqrt .5 mul -12 add}
\psplot[linecolor=blue,plotpoints=100]{0}{15}{x 2 mul 2 exp 1 add sqrt .6 mul 999 sqrt 0.8 mul add .5 mul -12 add}
\psline[arrowsize=.2 2]{->}(0,0)(0,10)
\psline[arrowsize=.2 2]{->}(0,0)(15,0)
\put(0.3,9.5){$t$}
\put(14.5,-0.8){$q_y$}
\put(-1,3.7){\red  $t_{0}$}
\put(-1.2,.9){\blue $t_{h_1}$}
\put(3.6,5){\red $q_0(t)$}
\put(3.4,2){\blue $q_{h}(t)$}
\put(1.3,8.5){\magenta  $\tilde t_0(q_{\max})=\tilde t_{h}(q_{\max})=t_{h_2}$}
\psline(0,1.5)(1.3,1.5)
\psline[arrowsize=.1 1]{->}(0,1.5)(.7,1.5)
\psline(0,2.5)(3.,2.5)
\psline[arrowsize=.1 1]{->}(0,2.5)(1.5,2.5)
\psline(0,3.5)(4.7,3.5)
\psline[arrowsize=.1 1]{->}(0,3.5)(3,3.5)
\psline(4.8,4.5)(6.4,4.5)
\psline[arrowsize=.1 1]{->}(5,4.5)(5.8,4.5)
\psline(7.5,5.5)(8.1,5.5)
\psline[linestyle=dashed,linecolor=magenta](0,7.8)(15,7.8)
\psline[linestyle=dashed,linecolor=magenta](0,3.8)(15,3.8)}
\end{picture}
  \caption{\label{acoustique23d:fig:5}Integration
    first over $t$ then over $q_y$ for the head wave}
\end{minipage}
\end{tabular}
\end{figure}

We thus have:
$$\tilde u^-_{sx,Ps}(x,0,z,s)=\int_{0}^{+\infty} u^-_{sx,Ps}(x,0,z,t)e^{-st}\,dt
$$
and we conclude by using the injectivity of the Laplace transform.
\section{Numerical illustration}
\label{Numerical}
\begin{figure}[htbp]
  \centering
\setlength{\unitlength}{.9cm}
  \begin{picture}(6,6.5)(2,0)
\psset{xunit=1.5cm}
\psframe(0,.25)(6,5.75)
\psline(0,3)(6,3)
\put(1.5,4.3){$\Omega^+$}
\put(1.5,2.1){$\Omega^-$}
\put(4.3,5.2){Source}
\psdot(3,4.5)
\psline{<->}(3,3)(3,4.45)
\put(3.5,4.){\unit{500}{\metre}}
\psline{<->}(5,3)(5,4.95)
\put(8.5,4.3){\unit{533}{\metre}}
\psline{<->}(5,3)(5,1.05)
\put(8.5,2){\unit{533}{\metre}}
\psline{<->}(3.05,4.5)(5,4.5)
\put(5.9,4.6){\unit{400}{\metre}}
\psdot(5,5)
\psdot(5,1)
\put(7.4,5.8){Receiver 1}
\put(7.4,0.7){Receiver 2}
  \end{picture}
  \caption{Configuration of the experiment}
\label{validation2d:fig:6}
\end{figure}
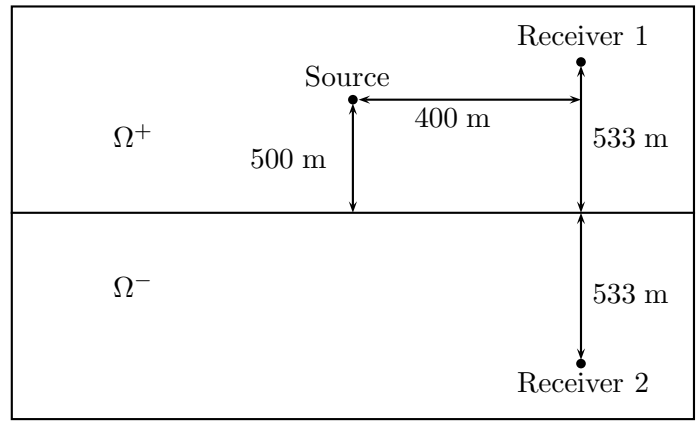

\begin{figure}[htbp]
      \begin{minipage}{.48\linewidth}
        \centerline{
          \setlength{\unitlength}{.9cm}
          \begin{picture}(8,6)
        \includegraphics[height=5cm]{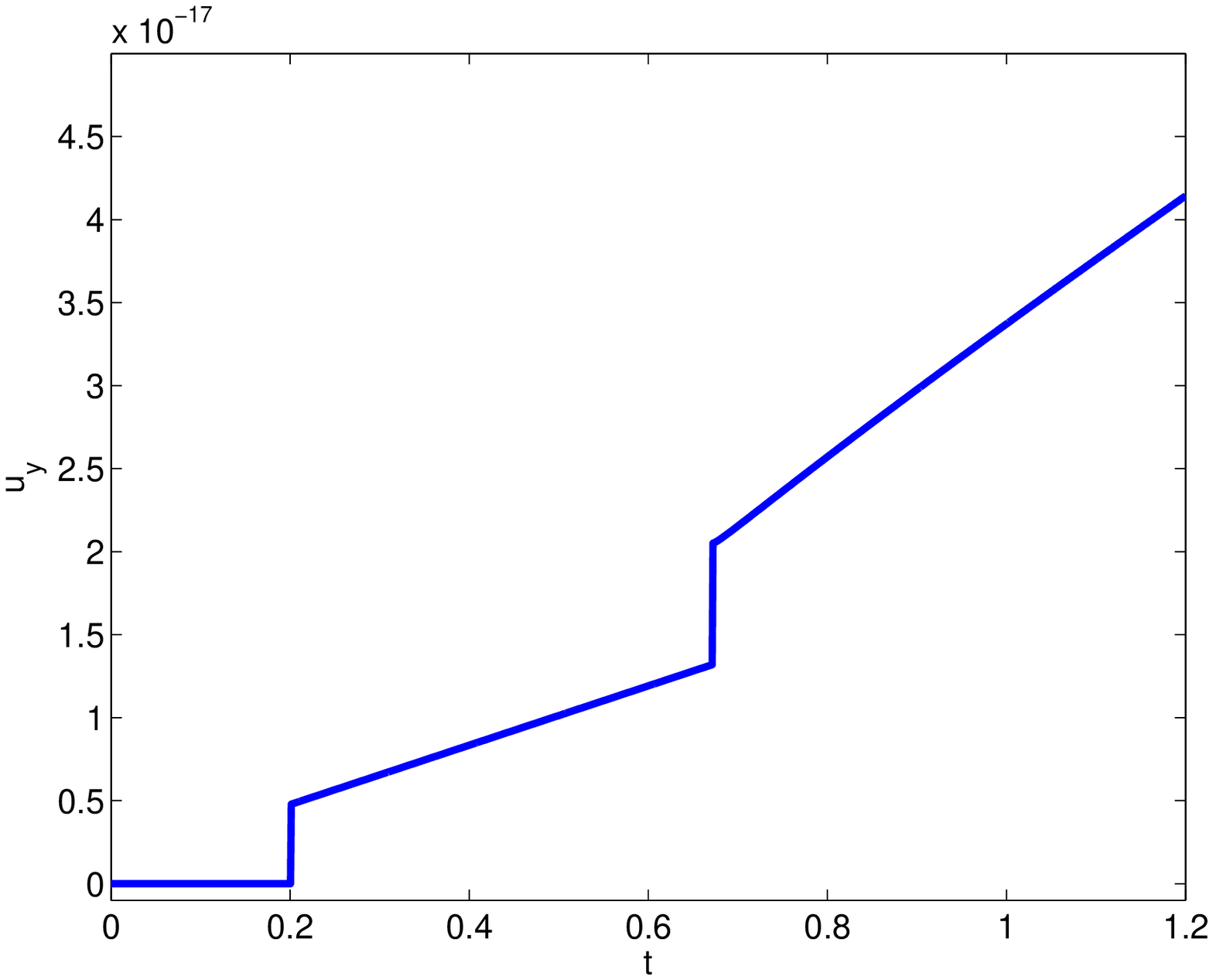}
\psline{<-}(-3.7,1.4)(-4.1,1.8)
\put(-5.4,2.2){Incident}
\psline{<-}(-2.,3)(-2.5,3.5)
\put(-3.7,4.1){Reflected}
\end{picture}
}  
    \end{minipage}
      \begin{minipage}{.48\linewidth}
      \centerline{
\setlength{\unitlength}{.9cm}
\begin{picture}(8,6)
        \includegraphics[height=5cm]{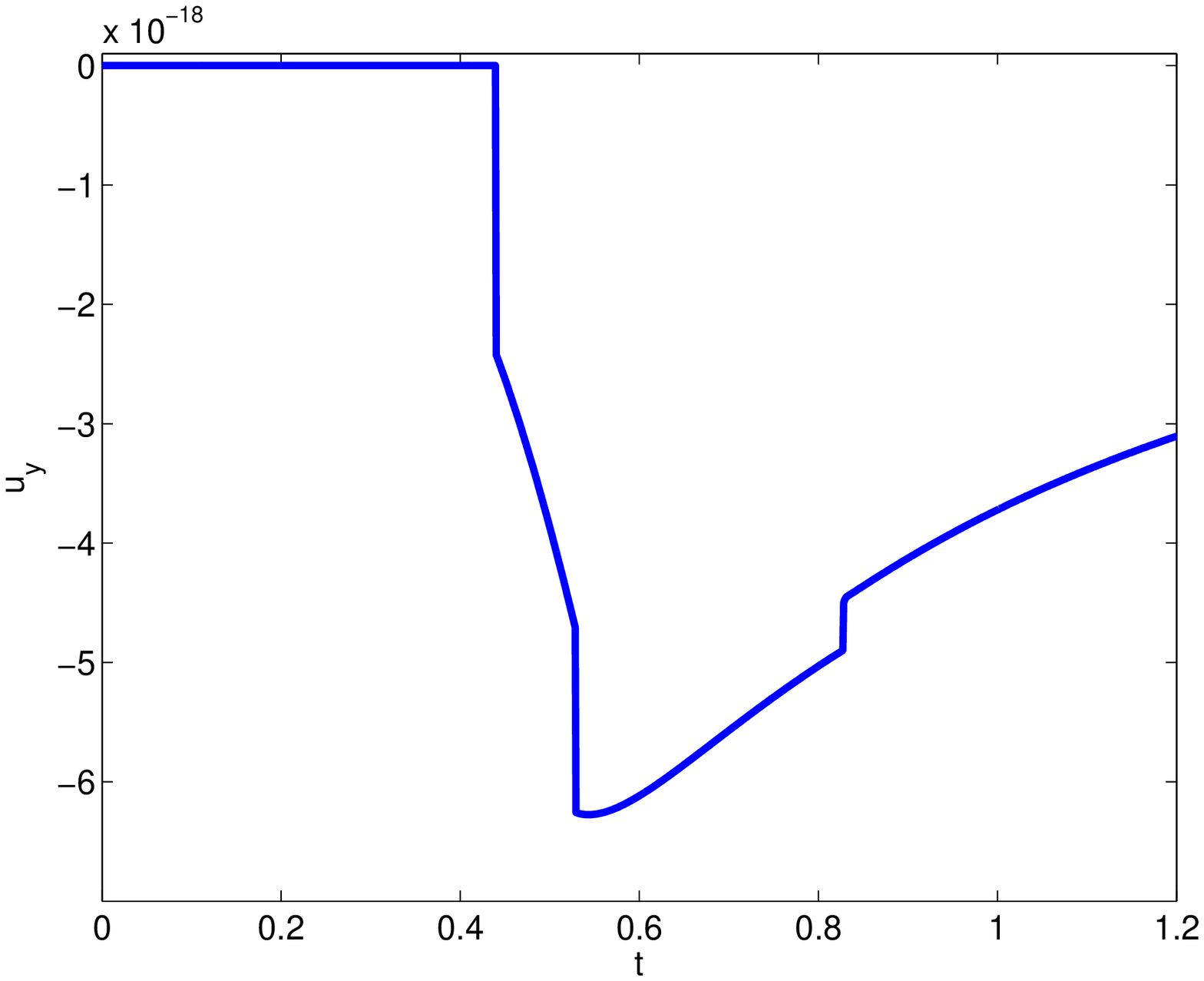}
\psline{<-}(-3.5,2.2)(-4.3,3)
\put(-5.1,3.5){$Pf$}
\psline{<-}(-2.5,1.1)(-1.45,1.1)
\put(-1.5,1.05){$S$}
\psline{<-}(-1.4,2.5)(-1.9,3.)
\put(-2.6,3.5){$Ps$}
\end{picture}
}
    \end{minipage}
    \caption{The $z$ component of the green
      function associated to the displacement at
      receiver 1 (left picture) and 2 (right
      picture).}
\label{acousporo:fig:2}
  \end{figure}

\begin{figure}[htbp]
      \begin{minipage}{.48\linewidth}
      \centerline{
          \setlength{\unitlength}{.9cm}
          \begin{picture}(8,6)
        \includegraphics[height=5cm]{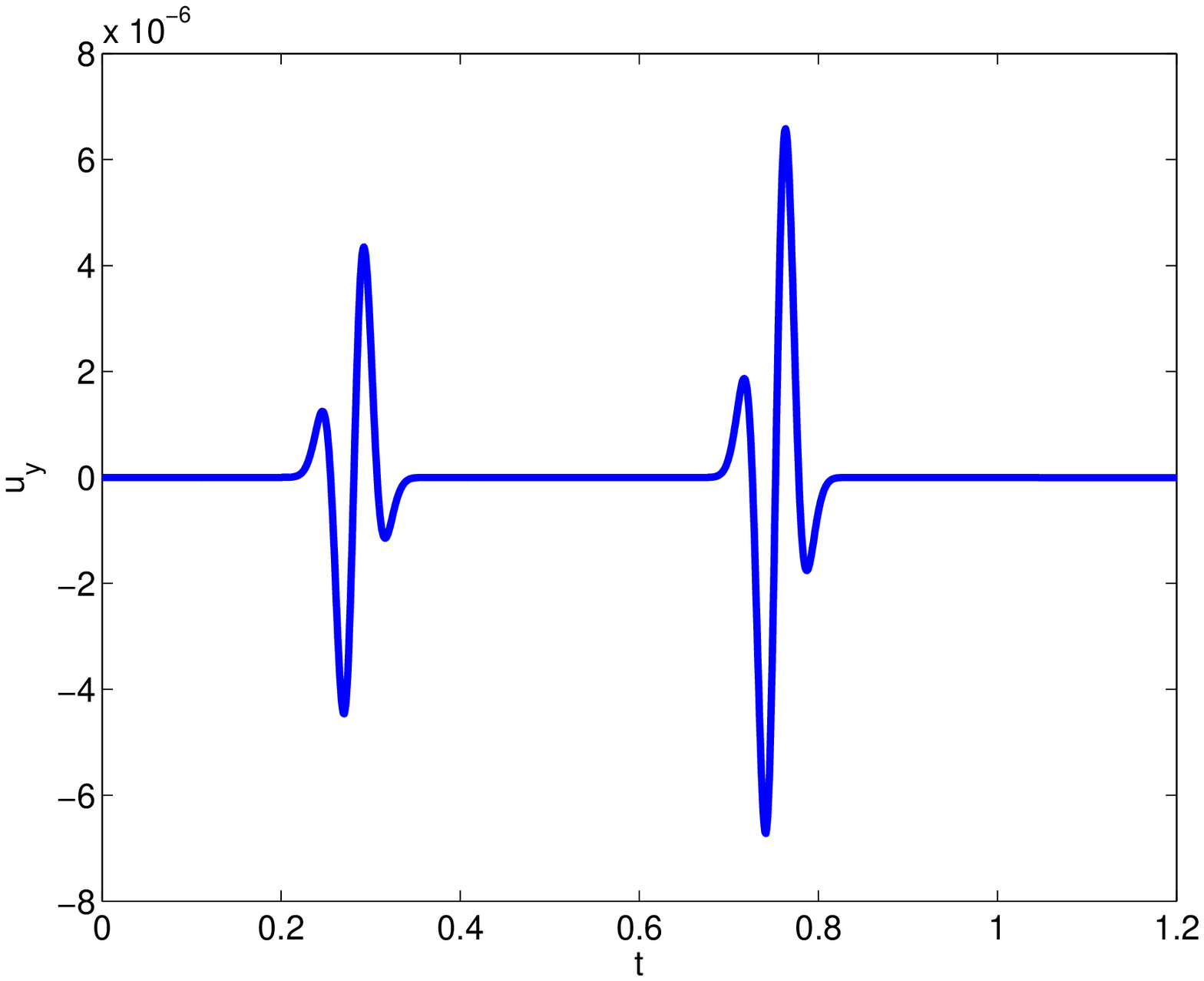}
\psline{<-}(-4.3,1.8)(-3.7,1.3)
\put(-4.5,1.){Incident}
\psline{<-}(-2.5,3.3)(-3.3,3.9)
\put(-4.5,4.5){Reflected}
\end{picture}
      }  
    \end{minipage} 
      \begin{minipage}{.48\linewidth}  
      \centerline{
          \setlength{\unitlength}{.9cm}
          \begin{picture}(8,6)
        \includegraphics[height=5cm]{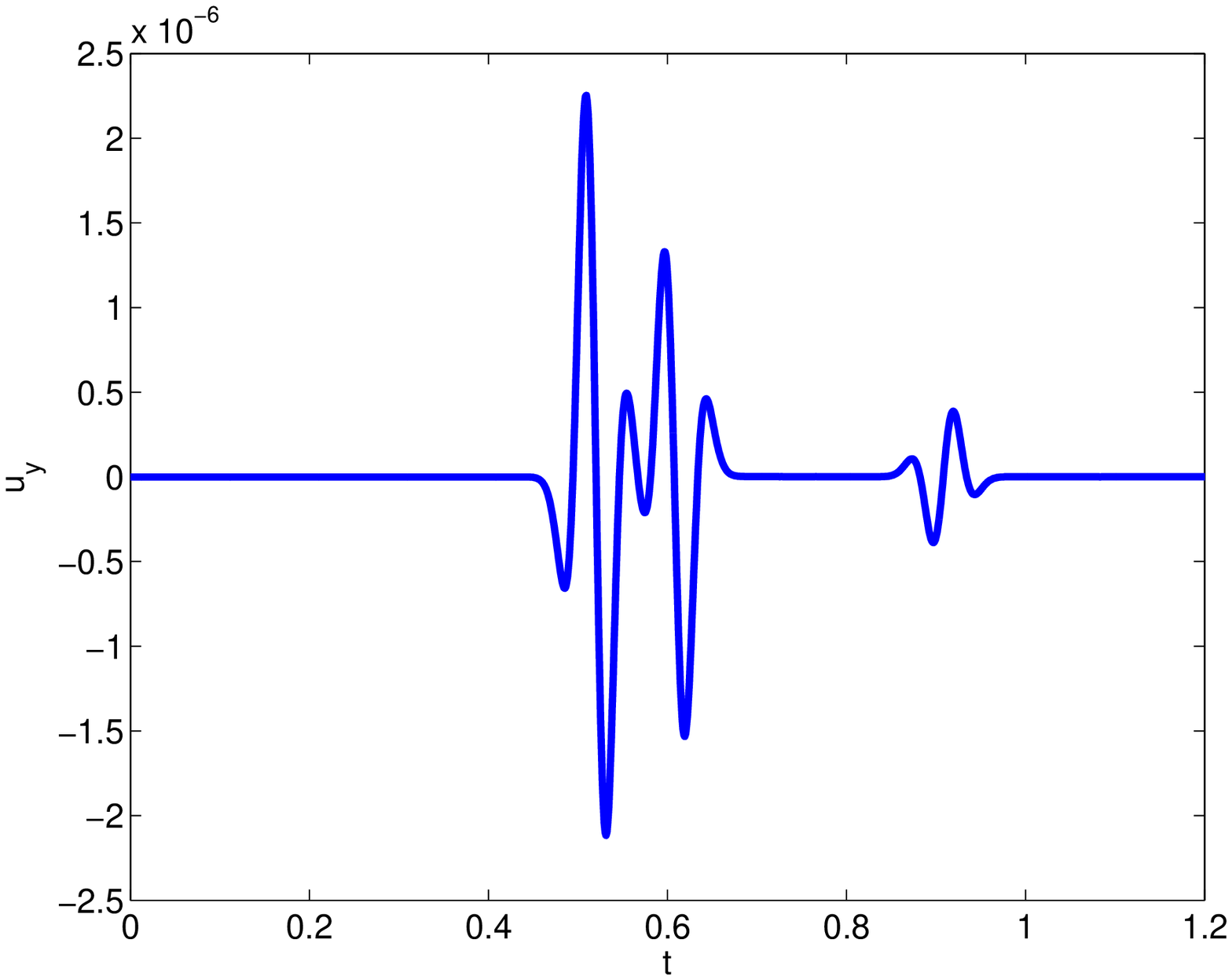}
\psline{<-}(-3.5,3.2)(-4.3,4)
\put(-5.05,4.55){$Pf$}
\psline{<-}(-2.8,1.5)(-1.9,1.5)
\put(-2.,1.5){$S$}
\psline{<-}(-1.6,3)(-2.1,3.5)
\put(-2.6,4){$Ps$}
 \end{picture}    }  
    \end{minipage}
\caption{The $z$ component of the displacement at
  receiver 1 (left picture) and 2 (right
  picture).}
\label{acousporo3d4:fig:1}
  \end{figure}

To illustrate our results, we have
computed the green function and the analytical solution to the following problem: we consider an
acoustic layer with a density
$\rho^+=\unit{1020}{\kilo\gram\per\cubic\metre}$
and a celerity
$V^+=\unit{1500}{\meter\per\second}$ on top of a
poroelastic layer whose characteristic
coefficients are:
\begin{itemize}
\item the solid density $\rho_s^-=\unit{2500}{\kilo\gram\per\cubic\metre}$;
\item the fluid density $\rho_f^-=\unit{1020}{\kilo\gram\per\cubic\metre}$;
\item the porosity $\phi^-=0.4$;
\item the tortuosity $a^-=2$;
\item the solid bulk modulus $K^-_s=\unit{16.0554}{\giga\pascal}$;
\item the fluid bulk modulus $K^-_f=\unit{2.295}{\giga\pascal}$;
\item the frame bulk modulus $K^-_b=\unit{10}{\giga\pascal}$;
\item the frame shear modulus $\mu^-=\unit{9.63342}{\giga\pascal}$;
\end{itemize}
so that the celerity of the waves in the
poroelastic medium are:
\begin{itemize}
\item for the fast P wave, $\Vpf^-=\unit{3677}{\meter\per\second}$
\item for the slow P wave, $\Vps^-=\unit{1060}{\meter\per\second}$
\item for the $\psi$ wave, $\Vs^-=\unit{2378}{\meter\per\second}.$
\end{itemize}
The source is located in the acoustic layer, at
\unit{500}{\meter} from the interface. It is a
point source in space and a fifth derivative of a
Gaussian of dominant frequency
$f_0=\unit{15}{\hertz}$:
$$f(t)=2\frac{\pi^2}{f_0^2}\left[3+12\frac{\pi^2}{f_0^2}\left(t-\frac{1}{f_0}\right)^2+4
  \frac{\pi^4}{f_0^4}\left(t-\frac{1}{f_0}\right)^4\right]e^{-\frac{\pi^2}{f_0^2}\left(t-\frac{1}{f_0}\right)^2}.$$
We compute the solution at two receivers, the
first one is in the acoustic layer, at
\unit{533}{\meter} from the interface; the first
one is in the poroelastic layer, at
\unit{533}{\meter} from the interface; both are
located on a vertical line at \unit{400}{\meter}
from the source (see
Fig.~\ref{validation2d:fig:6}).  To compute the
integrals over $q$ and the convolution with the
source function, we used a classical mid-point
quadrature formula.\\\\
We represent the
$z$ component of the green function associated to
the displacement from $t=0$ to
$t=\unit{1.2}{\second}$ on
Fig~\ref{acousporo:fig:2} and the displacement in
Fig.~\ref{acousporo3d4:fig:1}.  The left picture
represents the solution at receiver 1 while the
right picture represents the solution at receiver
2.  As all the types of waves are computed independently, it is easy to distinguish all of them, as it is indicated in the figures. 
solution.\\[10pt]

\section{Conclusion}
In this paper we have provided the complete solution
(reflected and transmitted wave) of the
propagation of wave in a stratified 3D medium
composed of an acoustic and a poroelastic layer. In a forthcoming paper  
we will extend the method to the propagation of waves in bilayered 
poroelastic medium in three dimensions.
\appendix
\section{Definition of the fictitious and real arrival times of the volume waves.}
We detail in this section the computation of the
fictitious and real arrival times of the transmitted $Ps$ wave at point
$(x,0,z)$. For a given $q_y\in\setR$, we  first determine 
  fastest path of the wave from the source to the point
  $(x,0,z)$, travelling at a velocity $\mat V^+(q_y)$ in the upper layer and at a velocity 
$\mat V_{Ps}^-(q_y)$ in the bottom layer: we search
  a point $\xi_0$ on the interface between the two
  media which minimizes the function
$$  t(\xi)=\frac {\sqrt{\xi^2+h^2}}{\mat V^+(q_y)}+\frac {\sqrt{(x-\xi)^2+z^2}}{\mat V_{Ps}^-(q_y)}$$
(see Fig.~\ref{de-rr:fig:1}). This leads us to find $\xi_0$ such that
\begin{equation}
  \label{de-rr:eq:1}
  t'(\xi_0)=\frac {\xi_0}{\mat V^+(q_y)\sqrt{\xi_0^2+h^2}}+\frac {\xi_0-x}
{\mat V_{Ps}^-(q_y)\sqrt{(x-\xi_0)^2+z^2}}=0.
\end{equation}
From a numerical point of view, the solution of
this equation is done by computing the roots of
the following fourth degree polynomial
$$ \begin{array}{ll}
&\dsp
\left(\frac{1}{{\mat V^+}^2(q_y)}-\frac{1}{{\mat V_{Ps}^-}^2(q_y)}\right)X^4+2x\left(\frac{1}{{\mat V_{Ps}^-}^2(q_y)}-\frac{1}{{\mat V^+}^2(q_y)}\right)X^3\\[18pt]
+&\dsp\left(
  \frac{x^2+z^2}{{\mat V^+}^2(q_y)}
-\frac{x^2+h^2}{{\mat V_{Ps}^-}^2(q_y)}\right)X^2+\frac{xh^2}{{\mat V_{Ps}^-}^2(q_y)}X+\frac{x^2h^2}{{\mat V_{Ps}^-}^2(q_y)},
\end{array}
$$
$\xi_0$ is thus the only real root of this polynomial located between 0 and $x$ which is also solution of~\eqref{de-rr:eq:1}.
Once $\xi_0$ is computed, we can define
$$\tilde t_{0}(q_y)=\frac {\sqrt{\xi_0^2+h^2}}{\mat V^+(q_y)}+\frac {\sqrt{(x-\xi_0)^2+z^2}}{\mat V_{Ps}^-(q_y)}\hbox{ and } t_0= \tilde t_{0}(0)$$
\setlength{\unitlength}{1cm}
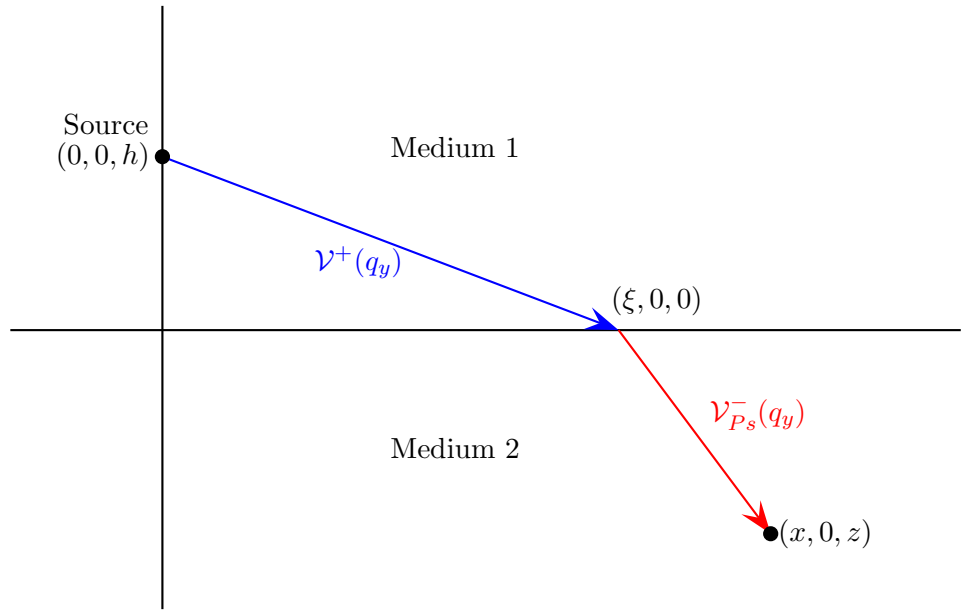
\begin{figure}[htbp]
\begin{picture}(10,10)
\put(2.7,6.8){Source}
\put(7.,6.5){Medium 1}
\put(7.,2.5){Medium 2}
\put(2.6,6.4){$(0,0,h)$}
\put(9.9,4.5){$(\xi,0,0)$}
\psline(2,4.2)(14.5,4.2)
\psline(4,0.5)(4,8.5)
\put(6,5.){$\blue \mat V^+(q_y)$}
\put(11.2,3.){$\red \mat V_{Ps}^-(q_y)$}
\psline[linecolor=blue,arrows=->,arrowsize=.2 4](4,6.5)(10,4.2)
\psline[linecolor=red,arrows=->,arrowsize=.2 4](10,4.2)(12,1.5)
\pscircle[fillstyle=solid,fillcolor=black](4,6.5){.1}
\pscircle[fillstyle=solid,fillcolor=black](12,1.5){.1}
\put(12.1,1.4){$(x,0,z)$}
\end{picture}
\caption{Path of the transmitted $i$ wave}
\label{de-rr:fig:1}
\end{figure}
Let us remark that 
\begin{proper}
  Since the
  fictitious velocities are smaller than the real
  one, the fictitious arrival times are greater than
  the real one. Moreover, since the fictitious
  velocities are even functions decreasing on $\setR^+$, $\tilde t_0$
  is an even function, increasing on
  $\setR^+$.
\end{proper}
\begin{coro}
  The function $\tilde t_{0}$ is
  bijective from $\setR^+$ to $\setR^+$.
\end{coro}
\bibliography{de-rr}
\bibliographystyle{plain}
\tableofcontents

\end{document}

%% file: ncacousporo.tex
\newcommand{\nc}{\newcommand}
\nc{\nl}{\newline}
\nc{\np}{\newpage}
\nc{\bequ}{\begin{equation}}
\nc{\eequ}{\end{equation}}
\nc{\bequs}{\begin{subequations}}
\nc{\eequs}{\end{subequations}}
\nc{\beqa}{\begin{eqnarray}}
\nc{\eeqa}{\end{eqnarray}}
\nc{\beqas}{\begin{eqnarray*}}
\nc{\eeqas}{\end{eqnarray*}}
\nc{\bcen}{\begin{center}}
\nc{\barr}{\begin{array}}
\nc{\earr}{\end{array}}
\nc{\ecen}{\end{center}}
\nc{\bite}{\begin{itemsize}}
\nc{\eite}{\end{itemsize}}
\nc{\benu}{\begin{enumerate}}
\nc{\eenu}{\end{enumerate}}
\nc{\btab}{\begin{tabular}}
\nc{\etab}{\end{tabular}}
\nc{\bfig}{\begin{figure}}
\nc{\efig}{\end{figure}}
\nc{\lea}{\left\{}
\nc{\rip}{\right.}
\nc{\itemb}{\item[$\bullet$]}
\nc{\ba}{\,/\,}
\nc{\dsp}{\displaystyle}
\nc{\mat}{\mathcal}
\nc{\be}{\beta}
\nc{\et}{\eta}
\nc{\si}{\sigma}
\nc{\Si}{\Sigma}
\nc{\la}{\lambda}
\nc{\La}{\Lambda}
\nc{\Gla}{\gr\La}
\nc{\ro}{\rho}
\nc{\om}{\omega}
\nc{\Om}{\Omega}
\nc{\al}{\alpha}
\nc{\ga}{\gamma}
\nc{\Ga}{\Gamma}
\nc{\To}{\tau}
\nc{\tz}{\tau_0}
\nc{\tu}{\tau_1}
\nc{\ts}{\tau_s}
\nc{\tp}{\tau_p}
\nc{\vare}{\varepsilon}
\nc{\vph}{\varphi}
\nc{\simu}{\si-\mu\frac{\partial u}{\partial x}}
\nc{\eps}{\epsilon}
\nc{\ups}{u_\eps}
\nc{\sps}{\si_\eps}
\nc{\p}{\partial}
\nc{\szs}{\scriptsize}
\nc{\RR}{\mbox{\rm I$\!$R}}
\nc{\CC}{\mbox{\rm I$\!$C}}
\nc{\NN}{\mbox{\rm I$\!$N}}
\nc{\KK}{\mbox{\rm I$\!$K}}
\nc{\LL}{\mbox{\rm I$\!$L}}
\nc{\MM}{\mbox{\rm I$\!$M}}
\nc{\II}{\mbox{\rm I$\!$I}}
\nc{\Ls}{\mat L^{sym}}
\nc{\un}{\underline}
\nc{\M}{\un{M}}
\nc{\GG}{\un{G}}
\nc{\HH}{\un{\un{H}}}
\nc{\XX}{\un{\un{X}}}
\nc{\XS}{\XX^{sym}}
\nc{\Xh}{\un{\un{X_h}}}
\nc{\Xhs}{\Xh^{sym}}
\nc{\Mh}{\un{M_h}}
\nc{\vv}{\un{\un{v}}}
\nc{\gr}{\boldsymbol}
\nc{\fr}{\frac{1}{2}}
\nc{\diva}{\mbox{\bf div}\,}
\nc{\tr}{\mbox{\bf Tr}}
\nc{\intt}{\int_{\scriptsize{Q_T}}}
\nc{\intr}{\int_{\RR}}
\nc{\Rd}{\RR^d}
\nc{\intn}{\int_{\scriptsize{\Rd}}}
\nc{\intd}{\intn}
\nc{\into}{\int_{\Om_\nu^t}}
\nc{\intg}{\int_{\Ga_\nu^t}}
\nc{\Lts}{L^2(\Rd,\Ls(\Rd))}
\nc{\Ld}{[L^2(\Rd)]^d}
\nc{\Hd}{[H^1(\Rd)]^d}
\nc{\blak}{\blacksquare}
\nc{\C}{\gr C}
\nc{\D}{\gr D}
\nc{\G}{\gr G}
\nc{\R}{\gr R}
\nc{\J}{\gr J}
\nc{\bn}{\big{\|}}
\nc{\gro}{\gr\ro}
\nc{\gaa}{\gr a}
\nc{\gm}{\gr m}
\nc{\gmt}{\gr{m_\To}}
\nc{\gb}{\gr b}
\nc{\difm}{\gr Z}
\nc{\Nabla}{\gr{\nabla}}
\nc{\U}{\mathbb{U}}
\nc{\V}{\mathbb{V}}
\nc{\W}{\mathbb{W}}
\nc{\Pe}{\mathbb{P}}
\nc{\lan}{\left<}
\nc{\ran}{\right>}
\nc{\ddw}{\ddot{\gr w}}
\nc{\ddwi}{\ddot{\gr w}_i}
\nc{\ddu}{\ddot{\gr u}_s}
\nc{\ddus}{\ddot{\gr u}_s}
\nc{\dduf}{\ddot{\gr u}_f}
\nc{\ddui}{\ddot{\gr u}_i}
\nc{\dw}{\dot{\gr w}}
\nc{\dwi}{\dot{\gr w}_i}
\nc{\du}{\dot{\gr u}_s}
\nc{\dui}{\dot{\gr u}_i}
\nc{\ddW}{\ddot{\gr W}}
\nc{\ddWi}{\ddot{\gr W}_i}
\nc{\ddU}{\ddot{\gr U}_s}
\nc{\ddUs}{\ddot{\gr U}_s}
\nc{\ddUf}{\ddot{\gr U}_f}
\nc{\ddUi}{\ddot{\gr U}_i}
\nc{\dW}{\dot{\gr W}}
\nc{\dWi}{\dot{\gr W}_i}
\nc{\dU}{\dot{\gr U}_s}
\nc{\duU}{\dot{\gr U}_i}
\nc{\ic}{{\rm i}\,}
\newtheorem{prop}{Proposition}[section]      %
\newtheorem{theo}{Theorem}[section]   %
\newtheorem{remark}{Remark}[section]          %
\newtheorem{lem}{Lemma}[section]             %
\newtheorem{defin}{Definition}[section]    %
\newtheorem{coro}{Corollaire}[section]       %
\newtheorem{hypo}{Assumption}[section]     %
\newtheorem{conc}{Conclusion}[section]       %
\newtheorem{proper}{Property}[section] 
\nc{\brem}{\begin{rem}}
\nc{\erem}{\end{rem}}
\nc{\btheo}{\begin{theo}}
\nc{\etheo}{\end{theo}}
\nc{\bprop}{\begin{prop}}
\nc{\eprop}{\end{prop}}
\nc{\blem}{\begin{lem}}
\nc{\elem}{\end{lem}}
\nc{\bdefi}{\begin{defin}}
\nc{\edefi}{\end{defin}}
\nc{\bcor}{\begin{coro}}
\nc{\ecor}{\end{coro}}
\nc{\bconc}{\begin{conc}}
\nc{\econc}{\end{conc}}
\nc{\bhyp}{\begin{hypo}}
\nc{\ehyp}{\end{hypo}}
\nc{\dem}{\proofname. }
\nc{\comm}{{\bf Commentaire }\\}
\newlength{\boxwidth}
\nc{\boxfigure}[3]{
                \setlength{\boxwidth}{#2}
                \addtolength{\boxwidth}{.1in}

                {\centering{
                \framebox[\boxwidth]{
                        \begin{minipage}{#2}
                        #3
                        \end{minipage}
                }}}}
\newlength{\fullboxwidth}
\newcommand{\rf}{{{\cal R}}}
\renewcommand{\ts}{{{\cal T}_{Ps}}}
\newcommand{\tf}{{{\cal T}_{Pf}}}
\newcommand{\tpsi}{{{\cal T}_{S}}}

\newcommand{\ka}{{\kappa^+}}
\newcommand{\kas}{{\kappa_{Ps}^-}}
\newcommand{\kaf}{{\kappa_{Pf}^-}}
\newcommand{\kapsi}{{\kappa_{S}^-}}

\newcommand{\setR}{\hbox{I}\!\hbox{R}}
 \newcommand{\setC}{{\mathord{\mathbb C}}}

\definecolor{pivert}{rgb}{0.2,0.8,0.2}

\newcommand{\Vpf}{V_{Pf}}
\newcommand{\Vps}{V_{Ps}}
\newcommand{\Vs}{V_S}

%% file: RR-6595.bbl
\begin{thebibliography}{10}

\bibitem{biot1}
M.~A. Biot.
\newblock Theory of propagation of elastic waves in a fluid-saturated porous
  solid. \rm{I}.~low-frequency range.
\newblock {\em J. Acoust. Soc. Am}, 28:168--178, 1956.

\bibitem{biot2}
M.~A. Biot.
\newblock Theory of propagation of elastic waves in a fluid-saturated porous
  solid. \rm{II}.~higher frequency range.
\newblock {\em J. Acoust. Soc. Am}, 28:179--191, 1956.

\bibitem{biot3}
M.~A. Biot.
\newblock Mechanics of deformation and acoustic propagation in porous media.
\newblock {\em J.~Appl. Phys.}, 33:1482--1498, 1962.

\bibitem{Cag}
L.~Cagniard.
\newblock {\em {Reflection and refraction of progressive seismic waves}}.
\newblock McGraw-Hill, 1962.

\bibitem{carcione}
J.~M. Carcione.
\newblock {\em Wave Fields in Real Media : Wave propagation in Anisotropic,
  Anelastic and Porous Media}.
\newblock Pergamon, 2001.

\bibitem{DH}
A.~T. de~Hoop.
\newblock The surface line source problem.
\newblock {\em Appl. Sci. Res. B}, 8:349--356, 1959.

\bibitem{diaz_th}
J.~Diaz.
\newblock {\em Approches analytiques et numériques de problèmes de transmission
  en propagation d'ondes en régime transitoire. Application au couplage
  fluide-structure et aux méthodes de couches parfaitement adaptées}.
\newblock PhD thesis, Université Paris 6, 2005.
\newblock in french.

\bibitem{RAP_DE6509}
J.~Diaz and A.~Ezziani.
\newblock Analytical solution for wave propagation in stratified
  acoustic/porous media. part \rm{I}: the 2\rm{D} case.
\newblock Technical Report 6509, INRIA, 2008.

\bibitem{RAP3}
J.~Diaz and A.~Ezziani.
\newblock Analytical solution for wave propagation in stratified poroelastic
  medium. part~\rm{I}: the 2\rm{D} case.
\newblock Technical Report 6591, INRIA, 2008.

\bibitem{Gar6}
J.~Diaz and A.~Ezziani.
\newblock Gar6more 2d.
\newblock \\ \url{http://www.spice-rtn.org/library/software/Gar6more2D}, 2008.

\bibitem{Gar63d}
J.~Diaz and A.~Ezziani.
\newblock Gar6more 3d.
\newblock \\ \url{http://www.spice-rtn.org/library/software/Gar6more3D}, 2008.

\bibitem{ezziani_th}
A.~Ezziani.
\newblock {\em Modélisation mathématique et numérique de la propagation d'ondes
  dans les milieux viscoélastiques et poroélastiques}.
\newblock PhD thesis, Université Paris 9, 2005.
\newblock in french.

\bibitem{FJ}
S.~Feng and D.~L. Johnson.
\newblock High-frequency acoustic properties of a fluid/porous solid interface.
  ii. the 2d reflection green's function.
\newblock {\em J. Acoust. Sec. Am.}, 74(3):915--924, 1983.

\bibitem{QG}
Q.~Grimal.
\newblock {\em Etude dans le domaine temporel de la propagation d'ondes
  elastiques en milieux stratifiés ; modèlisation de la reponse du thorax a un
  impacts.}
\newblock PhD thesis, Université Paris12-Val de Marne, 2003.
\newblock in french.

\bibitem{PG}
Y.~Pao and R.~Gajewski.
\newblock {\em The generalized ray theory and transient response of layered
  elastic solids}, volume~13 of {\em Physical Acoustics}, chapter~6, pages
  183--265.
\newblock 1977.

\bibitem{VDH}
J.~H. M.~T. van~der Hijden.
\newblock {\em {Propagation of transient elastic waves in stratified
  anisotropic media}}, volume~32 of {\em North Holland Series in Applied
  Mathematics and Mechanics}.
\newblock Elsevier Science Publishers, 1987.

\end{thebibliography}
